\documentclass{IEEEtran}  
\usepackage{algorithm,bbm}
\usepackage[noend]{algpseudocode}
\makeatletter
\def\BState{\State\hskip-\ALG@thistlm}
\makeatother
\usepackage{xcolor}
\usepackage{amssymb,latexsym,amsmath,array,graphicx,mathtools}
\usepackage[utf8]{inputenc}
\usepackage[T1]{fontenc}
\usepackage[ampersand]{easylist}
\usepackage{hyperref}
\usepackage{enumerate}
\usepackage{stfloats,dsfont}
\usepackage{mathrsfs}

  \newtheorem{innercustomlem}{Lemma}
\newenvironment{customlem}[1]
  {\renewcommand\theinnercustomlem{#1}\innercustomlem}
  {\endinnercustomlem}
  
    \newtheorem{innercustomprob}{Problem}
\newenvironment{customprob}[1]
  {\renewcommand\theinnercustomprob{#1}\innercustomprob}
  {\endinnercustomprob}
  
  \newtheorem{innercustomrem}{Remark}
\newenvironment{customrem}[1]
  {\renewcommand\theinnercustomrem{#1}\innercustomrem}
  {\endinnercustomrem}
  
  \newtheorem{innercustomthm}{Theorem}
\newenvironment{customthm}[1]
  {\renewcommand\theinnercustomthm{#1}\innercustomthm}
  {\endinnercustomthm}

  \newtheorem{innercustomassumption}{Assumption}
\newenvironment{customassumption}[1]
  {\renewcommand\theinnercustomassumption{#1}\innercustomassumption}
  {\endinnercustomassumption}

\newenvironment{proof}{{\em Proof:}}{\hfill$\fbox{}$}
\usepackage{mathtools}

\usepackage{scalefnt} 
\usepackage{blindtext}
\newcommand\scalemath[2]{\scalebox{#1}{\mbox{\ensuremath{\displaystyle #2}}}}
\usepackage[ampersand]{easylist}

\usepackage{mathrsfs}
\DeclareMathAlphabet{\mathpzc}{OT1}{pzc}{m}{it}
\usepackage{subcaption}

\title{\LARGE \bf
Stabilization of Underactuated Linear Coupled Reaction-Diffusion PDEs via Distributed or Boundary Actuation
}
\author{Constantinos Kitsos and Emilia Fridman
\thanks{This work was supported by the Israel Science Foundation (grant no. 673/19) and by the Chana
and Heinrich Manderman Chair at Tel Aviv University. The work of Constantinos Kitsos was supported by the Bloomfield International Postdoctoral Fellowship Fund}

\thanks{The authors are with the School of Electrical Engineering, Tel-Aviv University, Tel Aviv, Israel (emails: \{constantinos, emilia\}@tauex.tau.ac.il).}

}
\begin{document}
\maketitle
\begin{abstract}
This work concerns the exponential stabiliza6 tion of underactuated linear homogeneous systems of $m$parabolic partial differential equations (PDEs) in cascade (reaction–diffusion systems), where only the first state is controlled either internally or from the right boundary and in which the diffusion coefficients are distinct. For the distributed control case, a proportional-type stabilizing control
 is given explicitly. After applying modal decomposition, the
 stabilizing law is based on a transformation for the ordinary
 differential equations (ODE) system corresponding to the
 comparatively unstable modes into a target one, where the
 calculation of the stabilization law is independent of the
 arbitrarily large number of these modes. This is achieved
 by solving generalized Sylvester equations recursively. For
 the boundary control case, under appropriate sufficient
 conditions on the coupling matrix (reaction term), the proposed controller is dynamic. A dynamic extension technique via trigonometric change of variables that places the
 control internally is first performed. Then, modal decomposition is applied followed by a state transformation of the
 ODE system, which must be stabilized in order to be written
 in a form where a dynamic law can be established. For
 both distributed and boundary control systems, a constructive and scalable stabilization algorithm is proposed, as the
 choice of the controller gains is independent of the number
 of unstable modes and only relies on the stabilization of
 the reaction term. The present approach solves the problem of stabilization of underactuated systems when in the
 presence of distinct diffusion coefficients, the problem is
 not directly solvable, similarly to the scalar PDE case. 
\end{abstract}

\section{Introduction}
The  control of systems of coupled parabolic PDEs in which not all states are controlled (underactuated systems) has attracted much attention and has been posed as an open problem in \cite{Zuazua (2007)}. Lions \cite{Lions (1989)} introduced the study of controllability of cascade systems of parabolic PDEs. Such systems have been studied theoretically meanwhile, (see survey \cite{Ammar-Khodja et al (2011)}, which collects the plethora of recent studies concerning various notions of controllability of underactuated coupled systems). As far as the stabilization problem for such systems is concerned, some answers to the problem of internal stabilization of cascaded parabolic systems via distributed control placed on one equation and acting in a part of the domain were given in our preliminary conference version of this article \cite{Kitsos and Fridman (2022)}. To the best of authors' knowledge, the stabilization of similar problems involving systems of $m$ parabolic PDEs in cascade with one controlled state either internally or from the boundary and in the presence of distinct diffusion coefficients has not been solved yet.

The motivation behind the class of interconnected parabolic systems we study comes from various areas including chemistry, electrophysiology, genetics, and combustion. More particularly, biological predator-prey models or population and social dynamics phenomena are modeled by coupled parabolic systems (``reaction-diffusion systems" \cite{Berryman (1992),Britton (1986)}, see also  the work in \cite{Wang et al (2020)}  for applications to information diffusion in social media,        the work in \cite{Gallinato et al (2017)} for tumor growth models, and the work in \cite{Meile and Jones (2016)} for microbial processes). Also, in chemical processes \cite{Orlov and Dochain (2002)}, coupled temperature-concentration parabolic PDEs are used to describe the process dynamics. In such systems, the requirement of not controlling all states arises naturally. Feedback stabilization for scalar parabolic PDEs has been studied intensively \cite{Boskovic et al (2001),Coron and Trelat (2004),Liu (2003),Russel (1978)}. For the vector case, boundary stabilization of parabolic systems where all states are controlled has been achieved via backstepping transformation \cite{Vazquez and Krstic} and via modal decomposition in \cite{Katz et al (2021)}. When underactuation is assumed in coupled PDEs, boundary stabilization has been achieved for some classes of hyperbolic systems in \cite{Aamo (2013),Coron et al (2013),Di Meglio et al (2013a)} via backstepping. For these hyperbolic systems, boundary stabilization is possible when some dissipativity property is fulfilled on the boundaries. For the class of parabolic systems, when considering underactuations, strict assumptions are imposed on the internal dynamics \cite{Baccoli et al (2015)}. In the latter, boundary stabilization was achieved for only two coupled parabolic PDEs with boundary control of the first state when a minimum-phase assumption is met in addition to other restrictions on plant and controller parameters. In that work, the stabilization of a parabolic system of more than two equations with fewer inputs than the number of states and in the presence of distinct diffusion coefficients was posed as an open problem. 

Various studies have been devoted to the controllability of underactuated systems with internal controls \cite{Burgos and de Teresa (2010),Crepeau and Prieur (2008),Duprez and Lissy (2016),Guerrero (2007)}. These manifest several difficulties which become more complicated with the number of the states and with the number of distinct diffusion coefficients as a result of the notion of \emph{algebraic solvability} \cite{Steeves et al (2019)}. The problem of internal stabilization (via distributed control) of such systems runs deep (see for instance \cite{Barbu (2010),Coron (2007),Munteanu (2019)}, see also \cite{Christofides (1998)}). Similar complications arise in boundary controllability and stabilizability studies \cite{Ammar-Khodja et al (2011)}. It is also revealed that the problem of boundary controllability for the vector parabolic case is significantly more difficult than the internal one (distributed control in a part of the domain) see \cite{Fernandez-Cara et al (2010)}. It turns out that an underactuated boundary control system is \emph{null controllable} when the diffusion coefficients are identical (see \cite[Th. 6.1]{Ammar-Khodja et al (2011)}). However, for distinct diffusion coefficients and more than two coupled equations, the problem becomes more intricate. Some solutions to distributed observer design problems with fewer observations than the number of the states and with distinct diffusion coefficients, and which demonstrate some degree of duality with regard to stabilization problems,  have been given in \cite{Kitsos et al (2021a)}, while in \cite{Kitsos et al (2021b)} and in \cite[Ch. 3]{Kitsos (2020)}, the cases of three and $m$ linear non-homogeneous hyperbolic coupled PDEs were studied. In these works, appropriate infinite-dimensional state transformations solving operator Sylvester equations were introduced to deal with distinct elements on the diagonal of the coefficient of systems' differential operators requiring the use of higher-order spatial derivatives as measurements to yield Lyapunov stabilization of the observer error. {We further refer to \cite{Alabau-Boussouira (2013)} for coupled systems in cascade. }

In this work, taking a step beyond the controllability studies, we solve the stabilization problem. We consider a system of $m$ parabolic PDEs in cascade with distinct diffusion coefficients, where only the first equation is controlled, and we follow a modal decomposition approach. For internal control, we generalize methods mainly used for the scalar case (see \cite{Coron and Trelat (2004)} on direct Lyapunov method for state feedback, see also \cite{Barbu et al (2006)}) to the case of underactuated systems with one scalar controller. We assume that the number of internal inputs appearing in the first equation is equal to the number of unstable modes and that the resulting matrix that multiplies control inputs in the unstable modes is nonsingular. We then introduce a novel state transformation for ordinary differential equations (ODEs) with dimension equal to the number of coupled PDEs and written as a polynomial matrix in the slower eigenvalues of a related Stürm-Liouville problem, with order related to the number of distinct diffusion coefficients. The coefficients of this polynomial matrix are nilpotent matrices up to the identity matrix, which are subject to recursive generalized Sylvester equations and can be easily determined via a provided algorithm, while their values depend on the dynamics of the parabolic system. The stabilizing law simply consists in determining control gains stabilizing the reaction matrix and also in calculating our introduced state transformation, which depends on system dynamics. In this way, for any given system specification we provide a construction of unified and scalable control laws independently of the number of eigenvalues needed to be stabilized, which can be arbitrarily large. For boundary control, we follow an indirect approach in order to place the controls internally and obtain a dynamic control law ({a PI controller)}. This is relevant since the Hautus test might fail in the presence of distinct diffusion coefficients when trying to stabilize directly by use of static feedback. We assume that the control placed on the right boundary of the first state is written as a sum of control components. Inspired by the recent dynamic extension approach in \cite{Karafyllis (2021)} for the scalar case, we adapt similar transformation to our vector case. The system is first mapped into a new one where the control components and their time-derivatives are placed internally in the PDEs. We then apply modal decomposition followed by another transformation to the eigenspectrum in order to place control components in the first equation. In the next step, we are in a position to choose the dynamic control law. It turns out that for distributed control, we can achieve an arbitrarily fast decay rate while for boundary control, this is not the case.  This work solves the problem of internal and boundary stabilization of underactuated systems, for which backstepping approaches have not been proven to give solutions yet and at the same time it provides a scalable stabilization algorithm despite the presence of distinct diffusion coefficients. 

{Our contribution is summarized by the following points:
(1) Constructive methods for a stabilization problem of underactuated coupled PDE systems. These include a scalable algorithm for the determination of a novel transformation based on Sylvester equations and PI controllers for the boundary control case via an introduced trigonometric extension. (2)  Introduction of a Sylvester-equation approach. Sylvester equations are widely used in the context of finite-dimensional systems and our method proposes the extension of such approaches to the context of PDE systems. (3) {Sufficient conditions and a solution to the boundary stabilization} problem of underactuated parabolic systems. The problem that we have identified was previously characterized as open and no solution has yet been achieved via backstepping \cite{Baccoli et al (2015)}.
}

The rest of this article is organized as follows. The system and the description of the problem of stabilization are presented in Section \ref{sec:problem}. The internal stabilization approach is presented in Section \ref{sec:internal}, where Theorem \ref{theorem1} provides the main relevant result. Section \ref{sec:boundary} concerns the boundary stabilization problem, where Theorem \ref{proposition} provides its solution. In Section \ref{sec:simul},  we  apply  our  methodology  to unstable reaction-diffusion systems by applying either distributed or boundary control, and in Section \ref{conclusion} we provide some conclusions. 

\textit{Notation:} For a given $x \in {{\mathbb{R}}^{m}}$, $\left| x \right|$ denotes its usual Euclidean norm and for a matrix $Q\in {{\mathbb{R}}^{m\times m}}$, ${Q}^{\top}$ denotes its transpose, $\left| Q \right|:=\sup \left\{ \left| Qw \right|,\left| w \right|=1 \right\}$ is its induced norm, $\mathrm{Sym}(Q)=\frac{Q+Q^{\top}}{2}$ stands for its symmetric part and $\lambda_{\min}(Q)$, $\lambda_{\max}(Q)$ denote its minimum and maximum eigenvalue, respectively. By $\text{diag}\{ A_1, \ldots, A_m \}$ (or $\text{blkdiag}$) we denote the diagonal (or block diagonal) matrix with elements $A_1,\ldots,A_m$ scalars (or matrices). By $I_m$ we denote the identity matrix of dimension $m$. By $\otimes$ we denote the Kronecker product. For $f, g$ in $L^2\left( 0,L;\mathbb R^m\right)$, by $\left <f,g\right >$ we denote the inner product $\left <f,g \right >=\int_0^L f^\top (x)g(x) dx$ with induced norm $\Vert \cdot \Vert_{L^2\left( 0,L;\mathbb R^m\right)}$, where $L^2\left( 0,L;\mathbb R^m\right)$ denotes the space of equivalence classes of measurable functions $f:[0,L]\to\mathbb{R}^m$. By $\ell^2(\mathbb N;\mathbb R^m)$ we denote the Hilbert space of the square summable sequences $x=(x_n)_{n=1}^{+\infty}$. 
By $\mathds 1_{\omega}$ we denote the indicator function of the set $\omega$. By $\delta_{ij}$  we denote the Kronecker delta $\delta_{ij} = 1$, if $i = j$ and $\delta_{ij} = 0$, otherwise and 
$\lceil \cdot \rceil $ stands for the ceiling function.

\section{Problem Statement and Requirements}\label{sec:problem}
In this section, we present the underactuated system with its requirements and the stabilization problem.

Consider a system of $m$ coupled 1-D parabolic PDEs in a finite domain with control acting on the first state only, written as follows for $(t,x)$ in $[0,+\infty)\times(0,L)$:
\begin{subequations}\label{sys}
\begin{align}
    &z_t(t,x)=D z_{xx}(t,x)+Q z(t,x)+ \theta B \sum_{j=1}^Nb_j(x)u_j(t),\label{sys0}\\
   \label{BC} &\gamma_{11} z(t,0)+\gamma_{12} z_x(t,0)=0, \notag\\&\gamma_{21} z(t,L)+\gamma_{22} z_x(t,L)=\left (1-\theta \right )B\sum_{j=1}^N u_j(t),\\
   \label{IC} &z(0,x)=z^0(x).
\end{align}
\end{subequations}
 System's state is represented by $z=\begin{pmatrix}z_1& \ldots& z_m\end{pmatrix}^\top$. Diffusion matrix $D=\text{diag}\left\lbrace d_1,\ldots,d_m\right\rbrace$ consists of diffusion coefficients $d_1,\ldots,d_m>0$. The  coupling (reaction term) and control matrices $Q$ and $B$ are assumed to be of the form 
\begin{gather} \label{matrixQ}
Q= 
{
\begin{pmatrix}
q_{1,1}&&\cdots &&q_{1,m}\cr 
q_{2,1}&
 & &&\cr
0&&&&\vdots\cr
\vdots & \ddots &\ddots &\ddots \cr
&&&&\cr
0&\cdots & 0& q_{m,m-1}&q_{m,m}
\end{pmatrix}\nonumber
}, \quad B=\begin{pmatrix}
1\\0\\\vdots\\\vdots\\\\0
\end{pmatrix}.\notag
\end{gather}
Scalar control actions $u_1(t),\ldots, u_{N}(t)$ with $N$ to be determined later, act on the first equation and parameter $\theta$ taking values in $\{0,1\}$ determines whether the control is placed internally ($\theta=1$) or on the right boundary ($\theta=0$). Functions $b_1(\cdot),\cdots,b_{N}(\cdot)$ in $L^2(0,L)$ describe how the internal control actions are distributed in $[0,L]$ and are subject to some constraints given below. On the boundaries, we have $\gamma_{ij}\in \mathbb R$ satisfying $\gamma_{i1}^2+\gamma_{i2}^2\neq 0,\quad  i=1,2$ when $\theta=1$, whereas for $\theta=0$, we have the additional restriction that $\gamma_{12}=1-\gamma_{11}$ and $\gamma_{22}=1-\gamma_{21}$ with $\gamma_{11}, \gamma_{21} \in \{0,1\}$. Condition on $\gamma_{ij}$ restricts the type of boundary conditions to either Neumann or Dirichlet ones for the boundary control case, whereas for the distributed control case, we can have more general boundary conditions of mixed type. For boundary control, this stands as a sufficient condition for invertibility of a matrix which leads to stabilizability as it is revealed later.  
 Systems of the form \eqref{sys} can model for instance Turing instability \cite{Turing (1952)} and instability of slime mold amoebae aggregation  \cite{Keller and Segel (1970)}.
We make the following assumption:
\begin{customassumption}{1}\label{assumption}
The elements of the subdiagonal of $Q$ satisfy
\begin{gather}
q_{2,1}, q_{3,2}, \ldots, q_{m,m-1}\neq 0, \label{controllability}
\end{gather}
which stands as a controllability condition for the pair $(Q,B)$. 
\end{customassumption}
Before presenting the stabilization method, 
consider the following family of scalar Stürm-Liouville eigenvalue problems for each $i=1,\ldots,m:$
\begin{align}
\begin{aligned}    
    d_i  \varphi^{\prime \prime}(x)+\bar \lambda \varphi(x) =&0, \quad 0<x<L,\\
    \gamma_{11} \varphi(0)+\gamma_{12}\varphi^\prime(0)=&\gamma_{21} \varphi(L)+\gamma_{22}\varphi^\prime(L)=0, \label{SLprob}
    \end{aligned}
\end{align}
admitting a sequence of eigenvalues $\bar \lambda_{n,i}=d_i  \lambda_n$, where $\lambda_n$ are the eigenvalues of \eqref{SLprob} with $d_i=1$. This sequence of eigenvalues corresponds to a sequence of eigenfunctions $(\varphi_n)_{n=1}^{+\infty}$. The eigenvalues form an unbounded increasing and non-negative sequence while the eigenfunctions form a complete orthonormal system in $L^2(0,L)$. {Note here that although we can easily derive explicit formulas for eigenfunctions and eigenvalues of the above Stürm-Liouville problems in the case of Neumann or Dirichlet boundary conditions (when one of the pairs $(\gamma_{1,i},\gamma_{2,j}),i,j=1,2$ is zero), in the case of Robin boundary conditions we do not have such explicit formulas. However, we may get some estimates of the eigenvalues, see for instance \cite[Sec. 3.3.1]{Orlov (2020)}.} 

When the control is placed internally, we make the following assumption on shape functions $b_j(\cdot)$.
\begin{customassumption}{2}\label{assumption3}
Matrix 
\begin{align}
\mathcal B_{N\times N}:=&\begin{pmatrix}\mathcal B_1^\top\\\vdots\\\mathcal B_{N}^\top
\end{pmatrix} \label{mathcalB}
\end{align}
consisting of 
$
\mathcal B_n:=\begin{pmatrix}
b_{1,n}&\cdots&b_{N,n}
\end{pmatrix}^\top,$
which contain projections $b_{j,n}:=\int_0^L   b_j(x) \varphi_n(x) dx,$ 
is nonsingular. 
\end{customassumption}
Assumption \ref{assumption3} leads to a stabilizability property as it is shown in the following. Similar assumption appears in several works in the context of stabilization of scalar parabolic PDEs, see for instance \cite{Hagen and Mezic (2003)}.
The next assumption concerns only the case of boundary control, namely, when $\theta=0$.
\begin{customassumption}{3}\label{assumption2}
When $\theta=0$, there exist $\delta_0, k_Q>0$ such that
\begin{align}
\text{Sym}\left (Q\right )-D
  \text{diag}\{ k_Q, \lambda_1,\ldots,\lambda_1\}+\delta_0 I_m \preceq 0. \label{stabilizability2}
\end{align}
\end{customassumption}
\begin{customrem}{1} 
The abovementioned condition on matrix $Q$ restricts the class of unstable reaction terms, we are allowed us to consider when dealing with the boundary stabilization problem. An even stronger version of it appears in \cite{Baccoli et al (2015)} (see Condition 2 in Section V therein) standing as a sufficient condition to solve the boundary stabilization problem for underactuated systems of two coupled parabolic PDEs via backstepping method. The system considered there is similar to the one we consider here, but with the restriction that $m=2$ only, i.e., two equations. In that work, it is also assumed that for the reaction term taking the form of $Q$, $q_{22}$ is negative (a minimum phase assumption). In our case, a weaker condition of the form $q_{22}<d_2 \lambda_1$ would be sufficient to guarantee that Assumption \ref{assumption2} holds. Note that in the same work \cite{Baccoli et al (2015)}, it is concluded that for the case of $m>2$ coupled equations as in the system \eqref{sys} we consider here, the problem of boundary stabilization is open. Note also that we do not at all invoke Assumption \ref{assumption2} when performing internal stabilization via distributed control (see Section \ref{sec:internal} below).
\end{customrem}

The rapid stabilization problem we wish to solve in this work is stated as follows:
\begin{customprob}{1} \label{Problem} 
Suppose that Assumptions \ref{assumption}-\ref{assumption2} hold true. Then, determine stabilizing laws for the two following stabilization problems, the internal one ($\theta=1$) and the boundary one ($\theta=0$).
\begin{enumerate}[i]
  \item Case $\theta=1$: For any $\delta>0$, find $N \in \mathbb N$ and internal stabilizing laws $u_1(t), \ldots, u_N(t)$ such that for $z^0$ in $H^1\left (0,L;\mathbb R^m\right)$ satisfying compatibility conditions, solutions to \eqref{sys} satisfy \begin{align}
    \Vert z(t,\cdot)\Vert_{L^2\left( 0,L;\mathbb R^m\right)} \leq \ell e^{-\delta t} \Vert z^0(\cdot) \Vert_{L^2\left( 0,L;\mathbb R^m\right)}, \forall t \geq 0 \label{stabilityproblem}
\end{align}
with  $\ell>0$. 
  \item Case $\theta=0$: For some $\delta_0>0$ satisfying \eqref{stabilizability2}, find $N \in \mathbb N$ and boundary stabilizing dynamic laws for $u_1(t), \ldots, u_N(t)$ such that for $z^0$ in $H^2\left (0,L;\mathbb R^m\right)$ satisfying $\gamma_{11} z^0(0)+(1-\gamma_{11})\left (z^0\right )^\prime(0)=\gamma_{21}z^0(L)+(1-\gamma_{21})\left (z^0\right )^\prime(L)=0$, solutions to \eqref{sys} satisfy \eqref{stabilityproblem}  with $\delta$ substituted by $\delta_0$.
\end{enumerate}
\end{customprob}
Answers to both cases (i) and (ii) of Problem \ref{Problem} are given in sections \ref{sec:internal} and \ref{sec:boundary}.

\section{Internal Stabilization} \label{sec:internal}

In this section, we provide a solution to internal stabilization described by Problem \ref{Problem} (case $\theta=1$). We first apply modal decomposition. Then exploiting the fact that the eigenspectrum of our operator is partitioned into an unstable (or slow) part and a stable (or fast) one thanks to the countability and monotonicity of the eigenvalues, we focus on the stabilization of the comparatively unstable modes. To stabilize these modes, we introduce a state transformation aiming at a stabilization reduction from dimension $mN$ to dimension $m$. This transformation is given explicitly after solving a family of generalized Sylvester equations. Finally, by Lyapunov's direct method, we achieve to prove the stabilization result.  

\subsection{Modal Decomposition and Proportional-Type Controller}\label{subsec:modal}

We apply modal decomposition and we study the finite-dimensional system corresponding to the comparatively unstable modes. 

Each of the states of \eqref{sys} can be presented as 
\begin{align}
    z_i(t,\cdot)=\sum_{n=1}^\infty z_{i,n}(t) \varphi_n(\cdot), \quad i=1\ldots,m
\end{align}
with coefficients $z_{i,n}$ given by 
\begin{gather}
z_{i,n}=\left <z_i,\varphi_n\right >\label{zi}.
\end{gather}
Taking the time-derivative of \eqref{zi}, substituting dynamics \eqref{sys}, and integrating by parts, we get the following dynamics for $z_n=\begin{pmatrix} z_{1,n}&\cdots& z_{m,n}\end{pmatrix}^\top:$  
\begin{align}
    &\dot z_n(t)=\int_0^L z_t(t,x) \varphi_n(x)  dx\notag\\&=\left [D z_x(\cdot)\varphi_n(\cdot)-Dz(\cdot)  \varphi_n^{\prime} (\cdot)\right ]_0^L \notag\\&+\left (-\lambda_n D +Q\right )z_n(t)+B \sum_{j=1}^N u_j(t)  \int_0^L \varphi_n(x)  b_j(x) dx, \notag
\end{align}
which by virtue of homogeneous boundary conditions for $\varphi_n(x)$ and $z(t,x)$, is written as follows:
\begin{align}
     \dot z_n(t)=&\left (-\lambda_n D +Q\right )z_n(t) +B \sum_{j=1}^N b_{j,n} u_j(t). \label{zn_dot0}
\end{align}
Now, given a desired decay rate $\delta>0$, by taking into account the countability and monotonicity of eigenvalues of the elliptic operator,  we can always find a $N \in \mathbb N$  large enough such that 
\begin{align}
-\lambda_{N+1} D +\text{Sym}(Q)+\delta I_m <0.\label{lambdan11}
\end{align} 
By monotonicity of $\lambda_n$, the above implies that
\begin{align}
-\lambda_{n} D +\text{Sym}(Q)+\delta I_m <0, \quad \forall n\geq N+1. \label{lambdan22}
\end{align}  
Using the notation $Z=\text{col} \{z_1,\ldots,z_N\}\in \mathbb R^{mN}$, we obtain the following system corresponding to the finite-dimensional part of the eigenspectrum of the parabolic operator:
\begin{align}
\dot Z(t)=A Z(t)+\tilde B u(t),  \label{z^Nopen}
\end{align}
where $u(t):=\begin{pmatrix} u_1(t)&\cdots&u_N(t)\end{pmatrix}^\top\in \mathbb R^N$,
\begin{align} 
    A:=&\text{blkdiag}\{-\lambda_1 D+Q,\ldots,-\lambda_N D+Q \}, \label{Amatrix1}
\end{align}
and  $\tilde B\in \mathbb R^{mN\times N}$ is given by $$
\tilde B:= \text{col}\left\lbrace B \mathcal B_1^\top,\ldots, B\mathcal B_N^\top\right \rbrace=\left ( I_N \otimes B \right )\mathcal B_{N \times N}. $$
By invoking the Hautus lemma, it is easy to see that the pair $(A, \tilde B)$ is stabilizable under Assumption \ref{assumption3}.

We now seek for feedback controls of proportional type written as
\begin{align} \label{control0}
u_j(t)=&K_j Z(t),
\end{align}
where $K_{j}\in \mathbb R^{1\times mN}$ are controller gains to be found below.
Then, a direct stabilization approach of system \eqref{z^Nopen} would require to solve inequality
$
\text{Sym} \left (\tilde P(A+\tilde B K) \right)+\delta \tilde P\prec 0, $
where $\tilde P$ in $\mathbb R^{mN\times mN}$ is symmetric positive definite and $K:=\text{col}\left\lbrace K_1,\ldots,K_N\right \rbrace$. The above is written  in the design linear matrix inequalities (LMI) form 
\begin{align}
\text{Sym}\left ( A\tilde P_{-1} +\tilde B O\right )+\delta \tilde P_{-1} \prec 0, \label{LMIfull0}
\end{align} 
where we denote the unknowns $\tilde P_{-1}=\tilde P^{-1}$ and $O=K\tilde P_{-1}$. Then, the desired gain matrix is given by $K=O \tilde P_{-1}^{-1}.$
This LMI involves matrices of dimension $mN$. 

In this work, we aim at reducing the dimension of the stabilization from $mN$, which depends on the number $N$ of modes to be stabilized, to just the dimension $m$ of the coupled parabolic system, which is fixed. It turns out that this requirement of stabilization is not directly met as a consequence of the presence of distinct diffusion coefficients $d_i$. In fact, we seek for stabilizing actuations $u_j(t)$, whose calculation up to an inversion of matrix $\mathcal B_{N\times N}$ does not depend on the number of the modes $N$ but only on the number of system's equations $m$. Such property is important when dealing with large instabilities in the dynamics or when one would need to efficiently tune the decay rate. In other words, stabilization of \eqref{z^Nopen} should be based on the stabilization of an $m \times m$ matrix, namely, reaction matrix $Q$ and not on each of the diagonal elements of $A$, which can be arbitrarily many depending on the number of modes we need to stabilize at a given rate $\delta$. 

In the next subsection, we will show via examples why stabilization of \eqref{z^Nopen} is not directly implementable when diffusion coefficients are distinct.

At this point, let us   denote 
\begin{align}
K=\mathcal B_{N\times N}^{-1}\text{blkdiag}\{\bar K_1, \ldots, \bar K_N\}\label{Knproperty0}
\end{align}
with $\bar K_1, \ldots, \bar K_N\in \mathbb R^{1\times m}$ to be determined later.
Closing the loop by use of feedback control \eqref{control0} and after change of feedback control variables \eqref{Knproperty0}, $Z$ satisfies dynamics
\begin{align} \label{z^Nsys}
    \dot Z(t)=\left ( A+F\right )Z(t),
\end{align}
where $A$ is given by \eqref{Amatrix1} and
\begin{align} 
F:=& \text{blkdiag} \{ B \bar K_1,\ldots B \bar K_N\}.\label{A,barb}
\end{align}
This block diagonal form of closed-loop system \eqref{z^Nsys} will permit us to apply later a universal stabilization law for all blocks simultaneously as it  is shown in the following analysis.
\begin{customrem}{2} 
In more theoretical studies on controllability issues for such coupled parabolic systems (see survey \cite{Ammar-Khodja et al (2011)}), the control term is usually of the form $B \mathds 1_{\omega} U(t,x)$, with control $U(\cdot, \cdot)$ time and space-dependent and $\omega$ a given open subset of $[0,L]$. We could have alternatively posed the present problem in this setting, however, in practical applications, shape functions $b_j(x)$ are already given to be fixed (see, for instance, \cite{Christofides (1998)}) and we seek for stabilizing actuations depending exclusively on time as in the present analysis. In the first scenario, we would have chosen a proportional-type controller (see \cite{Barbu (2010)} (Chapter 2), see also \cite{Munteanu (2019)} (Chapter 9) for the signle PDE case) of the form
$U(t,x)=\sum_{j=1}^N \bar K_j \beta_j(x) \int_0^L z(t,x) \varphi_j(x) dx,$
where $\bar K_j \in \mathbb R^{1\times m}$ while functions $\beta_j(\cdot)$ are chosen to be written as
 $\beta_j(x)=\sum_{k=1}^N \beta_{jk} \varphi_k(x), \quad j=1,\ldots, N, $
with coefficients $\beta_{jk} \in \mathbb R$ satisfying
   $ \sum_{k=1}^N \beta_{jk} \int_0^L \mathds 1_{\omega} \varphi_k(x) \varphi_n(x) dx=\delta_{jn},$ for all $ j, n =1,\ldots,N.$ 
The previous equation, thanks to the linear independence of eigenfunctions $\varphi_n$, leads to a unique solution for unknown coefficients $\beta_{jk}$. This solution is represented as 
$\begin{pmatrix}  \beta_{11}&\cdots&\beta_{1N}\\\vdots&&\vdots\\\beta_{N1}&\cdots&\beta_{NN}\end{pmatrix}^\top= \left ( \int_0^L \mathds 1_{\omega}\varphi_i(x)\varphi_j(x) dx , {i,j=1,\ldots,N}\right)^{-1}$. Then, the finite-dimensional part of the eigenspectrum satisfies the same equations as in \eqref{z^Nsys} and we may follow a similar approach as the one presented below. 
\end{customrem}

\subsection{Problem of Stabilization of the Unstable Modes} \label{subsec:stabilization}

We present below some scenarios of stabilization of the finite-dimensional part of the eigenspectrum decomposition revealing its difficulty when diffusion matrix $D$ has distinct elements, i.e., when our system has distinct diffusion coefficients. 
Let us consider \eqref{z^Nsys}. To achieve exponential stability of this system with decay rate $\delta$, one would need to stabilize each of the components $-\lambda_n D+Q$ of the block diagonal matrix $A$ at this rate by choice of appropriate gains $\bar K_n$ as in \eqref{LMIfull0}-\eqref{Knproperty0}. However, this stabilization strategy would require stabilization of an $mN\times mN$ matrix, which is inefficient when $N$ becomes large.
In order to reduce the stabilization problem for all $N$ modes to just the stabilization of the coupling matrix $Q$, we need to follow an indirect strategy. Indeed, following a direct approach and trying to stabilize only matrix $Q$, one would choose gains
$ \bar K_n=K_Q, n=1,\ldots, N,$ where
$K_Q\in \mathbb R^{1\times m}$ is chosen such that a
Lyapunov matrix inequality of the form  
\begin{align}\label{lyapmatrix}
\text{Sym}\left ( P\left (Q+BK_Q\right )\right )+qP<0
\end{align}
is satisfied for $P\in \mathbb R^{m\times m}$
symmetric positive definite, which is nondiagonal, and some $q>0$. This is always possible due to the controllability of $(Q,B)$.
Then, to check asymptotic stability of system \eqref{z^Nsys}, choose Lyapunov function of the form
\begin{align} \label{lyapfunctinitial}
V_0(t)=\frac{1}{2}(z^N(t))^{\top} \bar P z^N(t)
\end{align}
with $ \bar P=I_N \otimes P$ consisting of $N$ diagonal blocks $P$.
Then, observe that $\text{Sym} \left (P\left ( -\lambda_n D+Q+BK_Q\right ) \right )$ appearing when taking the time-derivative of $V_0$ is of indefinite sign since $D$ and $P$ do not commute when $D$ has distinct diffusion coefficients and because $P$ is nondiagonal. This means that a stabilizing law chosen to stabilize $Q$ would not automatically lead to the stabilization of all the modes we need to stabilize at rate $\delta$. Note that this complication arising from the lack of a commutative property between the coefficient of the differential operator (the diffusion matrix $D$ here) and a Lyapunov matrix $P$ has been tackled in \cite{Kitsos (2020)}.
To understand how the number of distinct diffusion coefficients plays a role in the complexity of the problem, let us see the following examples.
\subsubsection*{Example 1 (one diffusion coefficient) }\label{exmpl1} Assume that all diffusion coefficients $d_i$ are identical, namely, $$d_1=d_2=\ldots=d_m.$$
 Then, the stabilization problem would be trivial. Indeed, the gains of the stabilization law \eqref{control0} via \eqref{Knproperty0} can be chosen  as $$\bar K_n=K_Q,$$ for all $n=1,\ldots,N$,
where $K_Q\in \mathbb R^{1\times m}$ is chosen such that Lyapunov inequality \eqref{lyapmatrix} is satisfied for $P$ symmetric positive definite and $q>0$ sufficiently large depending on the choice of the desired decay rate $\delta$. Then, by choice of Lyapunov function \eqref{lyapfunctinitial}, matrix $\text{Sym} \left (P\left ( -\lambda_n D+Q+B K_Q\right ) \right )=\text{Sym} \left (P\left ( -\lambda_n d_mI_m+Q+B K_Q\right ) \right )$ is always negative definite and the decay rate of system \eqref{z^Nsys} can attain value $\delta$ after appropriate choice of $q$, namely, $q \geq \delta-\lambda_1 d_m$. 

\subsubsection*{Example 2 (two diffusion coefficients) }\label{exmpl2} Let us now see the case where diffusion coefficients are identical up to the second one, namely, $$d_1\neq d_2=\ldots=d_m.$$
We choose gains $\bar K_n$ in \eqref{Knproperty0} given as in Example 1, but with an extra term, namely, 
\begin{align}
    \bar K_n=-G_n+K_Q, \quad \forall n\in \{1,\ldots,N\}, \label{K_nsigma=2}
\end{align}
where $G_n:=\lambda_n\left (
    d_2-d_1\right )B^\top$.
Again, the gain $K_Q\in \mathbb R^{1\times m}$ is chosen to satisfy Lyapunov inequality  \eqref{lyapmatrix} and then by choice of Lyapunov function \eqref{lyapfunctinitial}, system \eqref{z^Nsys} is stabilized at rate $\delta$. This is possible by noting that matrix $\text{Sym} \left (P\left ( -\lambda_n D+Q+B \bar K_n\right ) \right )$, which by \eqref{K_nsigma=2} is equal to $\text{Sym} \left (P\left ( -\lambda_n d_mI_m + Q+B K_Q\right ) \right )$, is negative definite by \eqref{lyapmatrix} and the decay rate of system \eqref{z^Nsys} can be equal to $\delta$ by appropriate choice of $q$, namely, $q \geq \delta-\lambda_1 d_m$. 
 
\subsubsection*{Example 3 (three diffusion coefficients) }\label{exmpl3}  We finally consider the case with $$d_2 \neq d_3$$ and let us consider for simplicity a $3 \times 3$ system ($m=3$).
Here, we might have $2$ or $3$ distinct diffusion coefficients and this stabilization problem turns to be more complicated than the previous ones. Indeed, to utilize the previously described Lyapunov stabilization for \eqref{z^Nsys}, in the absence of commutative property between $P$ and $D$, we perform a transformation of the form $y_n=T_n z_n, $ for $n=1,\ldots,N$ with 
\begin{align}
    T_n=I_3+\lambda_n \begin{pmatrix}
    0 & \kappa& 0\\0&0&0\\0&0&0
   \end{pmatrix}; \quad \kappa:=\frac{{d_3-d_2}}{q_{21}}. \label{T_n3}
\end{align}
Then, $Y=\text{col}\{ y_1,\ldots,y_N\}\in \mathbb R^{3N}$ satisfies
\begin{align}  
    \dot Y(t)=\left (\tilde A + \tilde F \right )Y(t),
\end{align}
where
\begin{align}
   \begin{aligned}
\scalemath{1}{\tilde A:=}&\scalemath{0.9}{\text{blkdiag}\{-\lambda_1 d_3 I_3+Q +BG_1,\ldots, -\lambda_N d_3I_3+Q +BG_N \}; }\\
G_n:=&\begin{pmatrix} \lambda_n \left ( d_3-d_1+\kappa q_{21}\right )\\ \lambda_n^2 \kappa\left ( d_1-d_2-\kappa q_{21}\right )+\lambda_n \kappa \left ( q_{22}-q_{11}\right )\\\lambda_n \kappa q_{23}\end{pmatrix}^\top \notag
\end{aligned}
\end{align}
and  $\tilde F:=\text{blkdiag}\{ B \bar K_1 T_n^{-1},\ldots, B \bar K_N T_n^{-1}\}$.
Then, the stabilizing gains are chosen to be of the form
\begin{align*}
\bar K_n=\left ( -G_n+K_Q \right ) T_n, \quad \forall n\in \{1,\ldots,N\},
\end{align*}
where the first term is needed to eliminate the undesired terms $B G_n$ and, as in the previous examples, $K_Q\in \mathbb R^{1\times m}$ is chosen to satisfy \eqref{lyapmatrix} with $q$ large enough, namely, $q \geq \delta-\lambda_1 d_3$. 

The abovementioned examples show that the problem of stabilization of an underactuated system is more intricate when diffusion coefficients are distinct, particularly when we have more than two distinct ones. In fact, index
\begin{align}\label{sigma}
\sigma:=&\min\left\lbrace i: d_i = d_j, \forall j=i,i+1,\ldots,m\right\rbrace
\end{align}
assigned to system \eqref{sys} is an indicator of the complexity of the stabilization problem. The larger the value of $\sigma$ is, the more complex is to determine the stabilization law. In our previous examples, for system in Example 1, $\sigma$ was equal to 1 (one diffusion), while in Example 2, $\sigma$ was equal to $2$. Example 3 with $\sigma=3$ provides us with intuition on an indirect strategy we should follow for systems with $\sigma>3$.
In the next section, considering all poisble values of $\sigma$, we provide a stabilization law by determining a state transformation similarly as in \eqref{T_n3} for $m=3$. 

\subsection{Stabilization Reduction and Main Result}\label{sec:transformation}

In this section, we aim at determining gains $K_n$ that lead to a closed-loop system, for which we can prove exponential stability. Our main goal is to reduce the problem of stabilization for the $mN \times mN$ system to a stabilization problem for system of dimension as large as $m$. We seek for a state transformation that transforms system \eqref{zn_dot0} into a target one where this type of control may be easily applied.

 Based on the previous reasoning, we present a target system which allows the derivation of the stabilizing law. Let us apply a transformation $y_n=T_n z_n$ to system \eqref{zn_dot0} with $T_n \in \mathbb R^{m\times m}$ an invertible polynomial matrix given by
\begin{align} 
T_n=\left\{ \begin{array}{ll} 
I_m+ \sum_{i=1}^{\bar \sigma} \lambda_n^i \bar T_i, \quad & 1\leq  n \leq N,\\I_m, \quad &n\geq N+1\end{array}\right., \label{T_n}
\end{align}
where $$\bar \sigma:=\min\{2\sigma-3,2m-4\}$$ with $\sigma$ given by \eqref{sigma} and $\lambda_n^i$ denoting the $i$-th power of $\lambda_n$. Note that $\mathcal T:=(T_n)_{n=1}^{+\infty}:\ell^2(\mathbb N;\mathbb R^m)\to \ell^2(\mathbb N;\mathbb R^m)$ is a bounded operator with bounded inverse. Matrices $\bar T_i\in \mathbb R^{m\times m}$ are assumed to be nilpotent of the form \eqref{matrixTi} shown at the bottom of the next page,
\begin{figure*}[b]
\begin{align}
\bar T_i = \begin{pmatrix}
0&\cdots&0&\kappa_{1,\lceil\frac{i}{2}\rceil+1}^{(i)} &\kappa_{1,\lceil\frac{i}{2}\rceil+2}^{(i)}&\cdots&\cdots&\kappa_{1,m}^{(i)} \cr
0&\cdots&0&0&\kappa_{2,\lceil\frac{i}{2}\rceil+2}^{(i)} &\cdots&\cdots&\kappa_{2,m}^{(i)}\cr
\vdots&&\ddots&&\ddots\cr
0&\cdots&&0&0&\kappa_{m-2-\lceil\frac{i}{2}\rceil,m-2}^{(i)}&\kappa_{m-2-\lceil\frac{i}{2}\rceil,m-1}^{(i)}&\kappa_{m-2-\lceil\frac{i}{2}\rceil,m}^{(i)}\cr
0&\cdots&&0&\cdots&0&\kappa_{m-1-\lceil\frac{i}{2}\rceil,m-1}^{(i)}&\kappa_{m-1-\lceil\frac{i}{2}\rceil,m}^{(i)}\cr
0&\cdots&&0&\cdots&0&0&0\cr
\vdots&&&\vdots&&\vdots&\vdots&\vdots\cr
0&\cdots&&0&\cdots&0&0&0
\end{pmatrix}, \quad i =1,\ldots,\bar \sigma \label{matrixTi}
\end{align}
\end{figure*}
where $\kappa^{(i)}_{j,k}$ are some constants to be determined explicitly in the following, which strictly depend on the dynamics of \eqref{sys} and not on $\lambda_n$. Note that superscripts ${(i)}$ appearing in $\kappa^{(i)}_{j,k}$ represent indices referring to each of the matrices $\bar T_i$, while their subscripts $(j,k)$ refer to their position in matrices $\bar T_i$.
By use of this transformation, we aim at obtaining a target system, which after injection of control \eqref{control0} and by use of \eqref{Knproperty0}, can be written in the closed-loop form
\begin{align}  \left\{ \begin{array}{ll}  \begin{aligned}
    \dot y_n(t)&=\left (-\lambda_n d_m I_m +Q +B G_n \right.\\&\left.+ B \bar K_nT_n^{-1}\right )y_n(t),  \qquad n \leq N,\\
    \dot y_n(t)&=\left (-\lambda_n D+Q\right )y_n(t) \\&+B \sum_{j=1}^N b_{j,n} K_j Z(t),  \quad  n\geq N+1
    \end{aligned}\end{array}\right.  \label{target_sys} 
\end{align}
with 
$G_n$ given by 
\begin{align}
 &G_n=-B^\top \left(\left ( Q-\lambda_n d_m I_m\right )\left( \sum_{i=1}^{\bar \sigma} \bar T_i \lambda_n^i\right )\right.  \notag\\&\left. +\left( \sum_{i=1}^{\bar \sigma} \bar T_i \lambda_n^i\right )\left ( \lambda_n D-Q\right )+\left ( D-d_mI_m\right )\lambda_n \right)T_n^{-1}.  \label{Gn}
\end{align}
Note that in target system \eqref{target_sys},  matrix $\lambda_n d_m I_m$ commutes with any matrix $P$ that we shall use to construct the Lyapunov functional later. The importance of this commutative property was revealed in the previous subsection. In addition, term $B G_n y_n(t)$ is undesired in the stabilization process but it will be canceled by use of the gains $\bar K_n$, similarly as it was done in examples 2 and 3 of the previous section. 

Let us now assume that $\bar T_i$ satisfy the following recursive generalized Sylvester equations
for all $ i\in \{1,\ldots,\bar \sigma\}$:
\begin{align} 
    \left(  I_m-BB^\top\right ) \left ( Q \bar T_i-\bar T_i Q+ \bar T_{i-1}\left ( D-d_m I_m \right )\right )=0, \label{Sylvester}
\end{align}
where we define $\bar T_0:=I_m$.  If \eqref{Sylvester} holds, then it is proven later that $T_n$ maps \eqref{zn_dot0}  to target system \eqref{target_sys}.  We obtain the following result on solutions to \eqref{Sylvester}:
\begin{customlem}{1} \label{lemma} 
If Condition \eqref{controllability} holds true, there exist matrices $\bar T_i$ of the form \eqref{matrixTi} shown at the bottom of the next page, satisfying generalized Sylvester equations \eqref{Sylvester}. Their components $\kappa^i_{j,k}$ are obtained explicitly by Algorithm \ref{algorithm} as follows.
\end{customlem}
\begin{proof}
We are in a position to directly determine solutions of this family of generalized Sylvester equations.  Thanks to the special structure of $\bar T_i$, we can apply an elimination procedure of each element of matrix $\left(  I_m-BB^\top\right ) \left ( Q \bar T_i-\bar T_i Q+ \bar T_{i-1}\left ( D-d_m I_m \right )\right )$ in a recursive manner. For each row, we start from elimination of its rightmost element and then we eliminate one by one all of its elements by moving one position to the left. The procedure initiates at the lowest row with nonzero elements and when all elements of the current row are eliminated leftwards, we recede to the rightmost element of one row before it and we continue the same procedure until all elements are eliminated. For each of these eliminations, we calculate an element $\kappa^i_{j,k}$ as the sole unknown in this entry, which is written as a linear combination of elements that have been already calculated in precedent eliminations. One can easily check that retrieving a sole unknown component $\kappa^i_{j,k}$ for each of these eliminations is a consequence of the  special structure of \eqref{matrixTi} and controllability condition \eqref{controllability}. More precisely, Algorithm \ref{algorithm} below describes in detail how to calculate each of the elements $\kappa^i_{j,k}$ of $\bar T_i$.
\begin{algorithm}
\caption{Calculation of transformation $T_n$}\label{algorithm}
\begin{algorithmic}[1]
\Procedure{Compute $\kappa_{j,k}^i,$ for all  $i \in \{1,\ldots,\bar \sigma\}, j\in \{1,\ldots,m-1-\lceil\frac{i}{2}\rceil\},k \in \{\lceil\frac{i}{2}\rceil+1,\ldots,m\}.$}{}
\State {$\bar T_0:=I_m$ and matrices $\bar T_i$ have the form \eqref{matrixTi} for $i\in \{1,\ldots,\bar \sigma\}.$ }
\State{$i=1$.}
\While{$i\leq \bar\sigma,$}
\State{$j=m-1-\lceil\frac{i}{2}\rceil$.}
\While{$j \geq 1$}
\State{$k=m$}
\While{$k\geq j+\lceil\frac{i}{2}\rceil$}
\State{Calculate $\kappa^i_{j,k}$ by eliminating element $(j,k)$ of matrix $\left(  I_m-BB^\top\right ) \left ( Q \bar T_i-\bar T_i Q+ \bar T_{i-1}\left ( D-d_m I_m \right )\right )$. In each step, all calculated $\kappa^i_{j,k}$ are written as linear combinations of $\kappa^i_{j,k}$ already calculated in previous steps.}
\State {$k\gets k-1$.}
\EndWhile
\State {$j\gets j-1$.}
\EndWhile
\State {$i\gets i+1$.}
\EndWhile
\EndProcedure
\end{algorithmic}
\end{algorithm}
By applying this algorithm, we achieve to calculate all constants $\kappa^i_{j,k}$ appearing in \eqref{matrixTi}. Indeed, one can see that eliminating each of the elements of $\left(  I_m-BB^\top\right ) \left ( Q \bar T_i-\bar T_i Q+ \bar T_{i-1}\left ( D-d_m I_m \right )\right )$ in the exact order the abovementioned algorithm suggests, we obtain a corresponding equation of the form
\begin{align}
    &q_{j+1,j}\kappa^i_{j,k}- \sum_{l=0}^{m-j-1-\frac{\lceil{i}\rceil}{2}} \kappa^i_{j+1, j+\frac{\lceil{i}\rceil}{2}+1+l}q_{j+\frac{\lceil{i}\rceil}{2}+1+l,k}\notag\\&+\sum_{l=0}^{m-j-1-\frac{\lceil{i}\rceil}{2}}q_{j+1,j+1+l}\kappa^i_{j+1+l,k}+\kappa^{i-1}_{j+1,k}(d_m-d_{k})=0.\notag
\end{align}
Each of the abovementioned equations has a unique solution for $\kappa^i_{j,k}$ by virtue of controllability condition \eqref{controllability}. We therefore directly obtain for all $i \in \{1,\ldots,\bar \sigma\}, j\in \{1,\ldots,m-1-\lceil\frac{i}{2}\rceil\},k \in \{\lceil\frac{i}{2}\rceil+1,\ldots,m\}$ the following formula:
\begin{align}
    &\kappa^i_{j,k}=\frac{1}{q_{j+1,j}}\left ( \sum_{l=0}^{m-j-1-\frac{\lceil{i}\rceil}{2}} \kappa^i_{j+1, j+\frac{\lceil{i}\rceil}{2}+1+l}q_{j+\frac{\lceil{i}\rceil}{2}+1+l,k}\right.\notag\\&\left.-\sum_{l=0}^{m-j-1-\frac{\lceil{i}\rceil}{2}}q_{j+1,j+1+l}\kappa^i_{j+1+l,k}+\kappa^{i-1}_{j+1,k}(d_m-d_{k})\right ), \label{kappaformula}
\end{align}
where we define $\kappa^0_{j+1,k}:=\delta_{j+1,k},$ for all $j\in \{1,\ldots,m-2\}, k\in \{1,\ldots,m\}$. The above result can be verified by invoking induction arguments.
\end{proof}

Now, we are in a position to establish the following result on internal stabilization ($\theta=1$).
\begin{customthm}{1}\label{theorem1} 
Consider parabolic system \eqref{sys} with distributed control ($\theta=1$) and initial condition $z(0,\cdot)=:z^0(\cdot)\in H^1\left (0,L;\mathbb R^m\right)$ satisfying compatibility conditions. Assume that both controllability condition on $(Q,B)$ in Assumption \ref{assumption} and on shape functions  $b_j(\cdot)$ in Assumption \ref{assumption3} hold true. Consider also polynomial matrices $T_n$ given by \eqref{T_n} with $\bar T_i$ solving \eqref{Sylvester}.  Given a decay rate $\delta>0$, let $N\in \mathbb N$ be subject to \eqref{lambdan11}. Assume that there exist $0\prec P_{-1}\in  \mathbb R^{m \times m}$ and $J\in \mathbb R^{1\times m}$ satisfying the following LMI:
\begin{align} \label{omega,psi}
 \text{Sym}\left (QP_{-1}+B J\right )+\left(\delta-\lambda_1 d_m \right ) P_{-1} \prec 0.
\end{align} 
Denote $K_Q=JP_{-1}^{-1}$. Let $\bar K_n\in \mathbb R^{1\times m}$ be given by
$\bar K_n=\left (-G_n+K_Q \right)T_n,$
for all $n=1,\ldots, N,$ where via expression \eqref{Gn} for $G_n$, the above is rewritten as
\begin{align}
\bar K_n=B^\top \Big ( \left ( Q-\lambda_n d_m I_m\right ) T_n+T_n\left ( \lambda_n D-Q\right ) \Big ) +K_Q T_n.\label{gains0}
\end{align}
Then, controller \eqref{control0} with gains $K_n$, defined by \eqref{Knproperty0}, exponentially stabilizes \eqref{sys} with a decay rate $\delta$, meaning that the solutions of the closed-loop system satisfy the following inequality: 
\begin{align}
    \Vert z(t,\cdot)\Vert_{L^2\left( 0,L;\mathbb R^m\right)} \leq \ell e^{-\delta t} \Vert z^0(\cdot) \Vert_{L^2\left( 0,L;\mathbb R^m\right)}, \forall t \geq 0 \label{stability0}
\end{align}
with  $\ell>0$.  Moreover, \eqref{omega,psi} is always feasible. 
\end{customthm} 
\begin{proof}
See Appendix \ref{appendix1}.
\end{proof}

This result illustrates that stabilization just requires the determination of a stabilizing gain $K_Q$ for coupling matrix $Q$ and also the calculation of a family of nilpotent matrices $\bar T_i \in \mathbb R^{m\times m}$, whose number depends on the number of distinct diffusion coefficients $d_i$ (represented by $\bar \sigma$) while their values only depend on system's dynamics. These matrices $\bar T_i$ are calculated easily by following the algorithm Lemma \ref{lemma} suggests. Note for instance that for $3\times 3$ systems, the sole matrix $\bar T_i$ needed has a single element [see \eqref{T_n3}]. This stabilization method is scalable up to the inversion of matrix $\mathcal B_{N\times N}$, given by \eqref{mathcalB}, meaning that after  stabilizing matrix $Q$, if we want to change the number $N$ of modes to stabilize, it is not required to stabilize a new matrix. Note also, that as in the description of part (i) of Problem \ref{Problem}, we achieve stabilization at any decay rate $\delta$.
\begin{customrem}{3} 
The finite-dimensional transformation $T_n$ \eqref{T_n} is directly related to an infinite-dimensional one firstly introduced in \cite{Kitsos (2020),Kitsos et al (2021a)} to solve an observer design problem corresponding to various classes of coupled PDEs. In these works, the corresponding transformation was a matrix operator with high-order differentiations in its domain and being a solution of a Sylvester operator equation.  Note that those works captured space-varying and nonlinear dynamics. Such cases, being more general than the ones here, required strong regularity assumptions and cannot be tackled by modal decomposition.  
\end{customrem}
\section{Boundary Stabilization} \label{sec:boundary}
In this section, we consider the boundary stabilization of \eqref{sys}, in order to give an answer to part (ii) of Problem \ref{Problem} (case $\theta=0$), recalling also that for this case we assumed that on the boundaries we have $\gamma_{12}=1-\gamma_{11}$ and $\gamma_{22}=1-\gamma_{21}$ with $\gamma_{11}, \gamma_{21} \in \{0,1\}$. The approach consists of a dynamic extension via trigonometric change of variables, then modal decomposition and finally, appropriately selecting the dynamic law that the control actuations satisfy. 

\subsection{Dynamic Extension}
In the following, we present the first step towards the boundary stabilization of \eqref{sys}. It consists of the application of a state transformation for dynamic extension followed by modal decomposition. We then perform another transformation to the eigenspectrum, in order to construct the dynamic control law. 
 
Below, we apply a state transformation in order to place the control internally. Such a type of transformation, but for the scalar PDE system, has been introduced in \cite{Karafyllis (2021)} leading to dynamic extension. We adapt this kind of transformation to our vector PDE system with one control. 
Let us first choose constants $\mu_j>0, j=1,\ldots, N$ with $\mu_j \neq \lambda_n$ for all $j=1,\ldots, N, n \in \mathbb N$. Let also $\psi_j(\cdot) \in C^2[0,L], j=1,\ldots,N$ be chosen to satisfy the following boundary-value problems,  which for given $\mu_j$ are uniquely solvable,
\begin{align}
\begin{aligned} \label{psidyn}
\psi_j^{\prime \prime}(x)+\mu_j \psi_j(x)=&0,  \quad 0<x<L,\\
\gamma_{11} \psi_j(0) +(1-\gamma_{11}) \psi^\prime_j(0)=&0, \\ \gamma_{21} \psi_j(L)+(1-\gamma_{21}) \psi_j^\prime(L)=&1,
\end{aligned}
\end{align}
recalling the restriction that $\gamma_{12}=1-\gamma_{11}$ and $\gamma_{22}=1-\gamma_{21}$ with $\gamma_{11}, \gamma_{21} \in \{0,1\}$ as we assumed for the case of boundary control ($\theta=1$).
It is convenient to choose $\mu_j$ such that 
\begin{align}
\label{muj}
\sqrt{\mu_j}=&\sqrt{\bar \mu_j}+2\mu_0 \frac{\pi}{L}, \quad j= 1,\ldots, N;\\
\sqrt{\bar \mu_j}:=&\left (1-\vert \gamma_{11}-\gamma_{21}\vert \right )\left ( j- \frac{1}{2}\right ) \frac{\pi}{L}\notag\\&+\vert \gamma_{11}-\gamma_{21}\vert j \frac{\pi}{L},
\end{align}
where $\mu_0\in \mathbb N$ is some parameter to be chosen later.
Then, the unique solutions to \eqref{psidyn} are given by 
\begin{align} \psi_j(x)=&(-1)^{j}\left ( \frac{1-\gamma_{21}}{\sqrt{ \mu_j}}+\gamma_{21}\right )\notag\\&\times\left ((1-\gamma_{11}) \cos \left ( \sqrt{ \mu_j}x\right )-\gamma_{11} \sin \left ( \sqrt{ \mu_j}x\right )\right ), \label{psi_j}
\end{align} 
for all $x\in [0,L]$.

Let us perform a transformation of the form
\begin{align}
w(t,x)=z(t,x) -B \sum_{j=1}^N \psi_j(x) u_j(t). \label{transfKar}
\end{align}
System \eqref{sys} is written in the new coordinates as
\begin{align} \label{wsystem}
&w_t(t,x)=D w_{xx}(t,x)+Q w(t,x)+QB \sum_{j=1}^N \psi_j(x)u_j(t)\notag \\&-d_1 B \sum_{j=1}^N \mu_j \psi_j(x) u_j(t)-B\sum_{j=1}^N \psi_j(x) \dot u_j(t)\notag \\
&\gamma_{11} w(t,0)+(1-\gamma_{11})w_x(t,0)=0, \notag \\ &\gamma_{21} w(t,L)+(1-\gamma_{21})w_x(t,L)=0.
\end{align}

In the next step, we aim at determining appropriate dynamic control law. To do this, we perform modal decomposition of \eqref{wsystem}. Let us consider ansatz continuously differentiable solutions $z(t,\cdot)$ in $L^2\left (0,L;\mathbb R^m \right )$ with $u_j(\cdot) \in C^1[0,+\infty)$ for all $j=1\ldots,N$. Existence of unique solutions to the closed-loop system and their regularity are proven later in Appendix \ref{app:thm2}. Thus, we are in a position to present each of the states of \eqref{wsystem} as 
\begin{align}
    w_i(t,\cdot)=\sum_{n=1}^\infty w_{i,n}(t) \varphi_n(\cdot), \quad i=1\ldots,m
\end{align}
with coefficients $w_{i,n}$ given by 
\begin{gather}
w_{i,n}=\left <w_i,\varphi_n\right >\label{wi}.
\end{gather}
Taking the time-derivative of \eqref{wi}, substituting dynamics \eqref{sys}, and integrating by parts, we obtain the following dynamics for $w_n=\left (w_{1,n},\ldots, w_{m,n}\right )^\top:$  
\begin{align}
    &\dot w_n(t)=\int_0^L w_t(t,x) \varphi_n(x)  dx\notag\\&=\left [D w_x(\cdot)\varphi_n(\cdot)-D w(\cdot)  \varphi_n^{\prime} (\cdot)\right ]_0^L\notag\\&+\left (-\lambda_n D +Q\right )w_n(t)+QB \sum_{j=1}^N \psi_{j,n} u_j(t)\notag\\&-d_1 B \sum_{j=1}^N \mu_j \psi_{j,n}u_j(t)-B\sum_{j=1}^N \psi_{j,n}\dot u_j(t) \notag
\end{align} 
\begin{align}
\text{with } \psi_{j,n}:=\left <\psi_j, \varphi_n\right >, \label{psijn0}
\end{align}
 which by virtue of boundary conditions for $\varphi_n(x)$ and $w(t,x)$, is written as follows:
\begin{align}
     \dot w_n(t)=&\left (-\lambda_n D +Q\right )w_n(t) +QB \sum_{j=1}^N \psi_{j,n} u_j(t)\notag\\&-d_1 B \sum_{j=1}^N \mu_j \psi_{j,n}u_j(t)-B\sum_{j=1}^N \psi_{j,n}\dot u_j(t). \label{wn_dot}
\end{align}
{We notice here that it is not yet straightforward to determine appropriate dynamic control law. This results from the fact that one would need to first cancel terms $u_j$ from each of the equations of system \eqref{wn_dot} by choice of the dynamics $\dot u_j$. This is not yet possible since in the second equation of \eqref{wn_dot}, although $u_j$ appear due to the term $QB \sum_{j=1}^N \psi_{j,n} u_j(t)$, their time-derivatives $\dot u_j$ do not appear at all. Therefore, let us perform a transformation  of the form }
\begin{align}
\bar w_n=&w_n-\begin{pmatrix}
0\\\Psi_n^\top r_1(t)\\\vdots\\\Psi_n^\top r_{m-1}(t)
\end{pmatrix} \label{transfbarw}
\end{align}
for all $n \in \mathbb N$, where $\Psi_n:=\begin{pmatrix}
\psi_{1,n}&\cdots&\psi_{N,n}\end{pmatrix}^\top$ and $r_i \in C^1\left ([0,+\infty);\mathbb R^N\right ),  i =1,\ldots,m-1$ are subject to appropriate dynamics to be determined later.
Then, system \eqref{wn_dot} is written in the new coordinates as
\begin{align}
\dot{\bar w}_n(t)=&\left (-\lambda_n D +Q \right )\bar w_n(t)+\left (-\lambda_n D +Q \right ) \begin{pmatrix}
0\\\Psi_n^\top r_1(t)\\\vdots\\\Psi_n^\top r_{m-1}(t)
\end{pmatrix} \notag\\&+QB \sum_{j=1}^N \psi_{j,n}u_j(t)-d_1 B\sum_{j=1}^N \mu_j \psi_{j,n}u_j(t)\notag\\&-\begin{pmatrix}
\Psi_n^\top \dot u(t)\\\Psi_n^\top \dot r_1(t)\\\vdots\\\Psi_n^\top \dot r_{m-1}(t)
\end{pmatrix} . \label{barwn0}
\end{align}
Before choosing the dynamic control law, we are in a position to show the following result. Let us first define 
\begin{align} \label{Psi}
\Psi_{N\times N}:=\begin{pmatrix}
\Psi_1^\top\\
\vdots\\
\Psi_N^\top
\end{pmatrix}.
\end{align}
\begin{customlem}{2} \label{Lemma_invertible} 
Assume that $\psi_j$ are of the form \eqref{psi_j} for all $j\in \{1,\ldots, N\}$. Then, $\Psi_{N\times N}$ is invertible for all $N \in \mathbb N$ and its inverse,  denoted by $\Psi_{N\times N}^{-1}=\left ( \chi_{i,k}; {i,k=1,\ldots,N}\right )$, has elements 
\begin{align}
\label{invPsi} &\chi_{i,k}=\frac{\sqrt{2} (-1)^k}{2L}\left ( \vert \gamma_{11}-\gamma_{21} \vert \left(\frac{\gamma_{21}}{\sqrt{\lambda_k}}-\frac{\gamma_{11}}{\sqrt{ \mu_i}}\right )\right.\notag\\&\left.+ \frac{1-\vert \gamma_{11}-\gamma_{21} \vert}{\sqrt{\lambda_k}}\left ( {\gamma_{11}}+\left (1-\gamma_{11}\right ){\sqrt{  \mu_i}}\right )\right )\bar \chi_{i,k};\\
&\bar{\chi}_{i,k}:=\frac{\displaystyle\prod_{\ell=1}^N \left (\lambda_k- \mu_{\ell} \right )\left (\lambda_{\ell}- \mu_{i}\right )}{\left ( \mu_i-\lambda_k\right ) \displaystyle\prod_{\substack{\ell=1 \\ \ell\neq k}}^N\left ( \lambda_{\ell}-\lambda_k\right ) \displaystyle\prod_{\substack{\ell=1 \\ \ell\neq i}}^N\left (  \mu_i- \mu_{\ell}\right ) },\notag
\end{align}
for $i, k =1,\ldots,N.$
\end{customlem} 
\begin{proof}
Recalling that in this section we have assumed $\gamma_{11}, \gamma_{21}\in \{0,1\}$, we can easily see that eigenvalues and eigenfunctions of the Stürm-Liouville problem  are of the form
\begin{align}\label{varphi}
\lambda_n=& \left (1-\vert \gamma_{11}-\gamma_{21} \vert \right ) \frac{n^2 \pi^2}{L^2}+\vert \gamma_{11}-\gamma_{21} \vert \left (n-\frac{1}{2}\right )^2\frac{\pi^2}{L^2}, \notag \\ \varphi_n(x)=&\sqrt{2} \left ((1-\gamma_{11})\cos\left ( \sqrt{\lambda_n} x\right )+\gamma_{11}\sin \left ( \sqrt{\lambda_n} x\right )\right ). 
\end{align}
We next calculate $\psi_{j,n}$ for all $j, n=1,\ldots, N$ by \eqref{psijn0} as
\begin{align}
&\psi_{j,n}= \frac{(-1)^n L \sqrt{2}}{ \mu_j-\lambda_n} \left( \vert \gamma_{11} -\gamma_{21} \vert\left ( \gamma_{21}\sqrt{\lambda_n}-\gamma_{11}\sqrt{ \mu_j}\right) \notag\right.\\&\left.+\left(1-\vert \gamma_{11}-\gamma_{21}\vert \right)\left (\gamma_{11} \sqrt{\lambda_n}+(1-\gamma_{11})\sqrt{\frac{{\lambda_n}}{{ \mu_j}}} \right )\right).
\end{align} 
Then, $\Psi_{N\times N}^\top=\left ( \psi_{j,n}; {j,n=1\ldots,N}\right )$ is written as 
\begin{align}
&\Psi_{N\times N}^\top=L\sqrt{2}\vert \gamma_{11}-\gamma_{21}\vert \left (\gamma_{21} C \text{diag}\left\lbrace (-1)^n\sqrt{\lambda_n}\right\rbrace_{n=1}^N \right .\notag\\&\left.- \gamma_{11} \text{diag}\{\sqrt{ \mu_j}\}_{j=1}^N C \text{diag}\{(-1)^n \}_{n=1}^N  \right. \bigg )\notag\\&+L\sqrt{2}\left (1-\vert \gamma_{11}-\gamma_{21} \vert \right) \left. \bigg (\gamma_{11}  C \text{diag}\{ (-1)^n\sqrt{\lambda_n}\}_{n=1}^N \right.\notag\\&\left.+(1-\gamma_{11})\text{diag}\left\lbrace \frac{1}{\sqrt{  \mu_j}}\right\rbrace_{j=1}^NC \text{diag}\{ (-1)^n\sqrt{\lambda_n}\}_{n=1}^N  \right ),
\end{align}
where $C:=\left ( \frac{1}{ \mu_j-\lambda_n}; {j,n=1,\ldots,N}\right )$ is a Hilbert-type matrix. We now invoke result in \cite[Lemma 2.1]{Trench and Scheinok (1966)}, which shows that $C$ is invertible with inverse explicitly calculated in (2.1) therein. Using this result, it is easily verifiable that the elements of the inverse $\Psi_{N\times N}^{-1}$ are of the form \eqref{invPsi}.
\end{proof}
For more general boundary  conditions like the Robin ones, an analogous result as in the one of this lemma would be harder to achieve.

Let us now denote $u(t):=\text{col}\left \lbrace u_1(t),\ldots,u_N(t)\right \rbrace \in \mathbb R^N, X(t):=\text{col}\{u(t),r_1(t),\ldots,r_{m-1}(t)\} \in \mathbb R^{mN}$. By virtue of Lemma \ref{Lemma_invertible}, we are in a position to construct dynamic control law of the form
\begin{align}
\dot X(t)=H X(t) -\left (B \otimes I_N \right )v(t) \label{dynlaw0}
\end{align}
with 
\begin{align}
H:=&-H_0+Q\otimes I_N;
 \label{Gamma}\\
H_0:=\text{blkdiag}&\left\lbrace d_1 M, d_2 \Psi_{N\times N}^{-1}\Lambda \Psi_{N\times N},\ldots,\right.\\&\left.d_m\Psi_{N\times N}^{-1}\Lambda \Psi_{N\times N}\right\rbrace, \notag\\
M:=&\text{diag}\{\mu_1,\ldots,\mu_N\}, \notag \\ \Lambda:=&\text{diag}\{\lambda_1,\ldots,\lambda_N\}.\notag
\end{align}
Matrix $\Psi_{N\times N}$ is given by \eqref{Psi} and $v(t):=\text{col}\{v_1(t),\ldots,v_N(t)\} \in \mathbb R^N$ is  a control input to be chosen appropriately later.
 Injecting the abovementioned dynamical law in system \eqref{barwn0}, we get
\begin{align}
\dot{\bar w}_n(t)=\left (-\lambda_n D +Q \right )\bar w_n(t)+B\sum_{j=1}^N \psi_{j,n} v_j(t),
\label{barwn1}
\end{align}
for all $n\in  \mathbb N$.

System \eqref{barwn1} is written in a form resembling to the one that would be derived after applying modal decomposition for a system with internal actuations $v_1,\ldots, v_N$ multiplied by shape functions $\psi_1,\ldots,\psi_N$, which are placed on the first equation only (see \eqref{zn_dot0} and the analysis of the previous section on internal stabilization). 

At this point, let us choose $k_Q$ and $\delta_0$ satisfying Assumption \ref{assumption2}. For this decay rate $\delta_0$, at which we able to stabilize \eqref{barwn1}, we can always find a $N \in \mathbb N$  large enough in such a way that 
\begin{align}
-\lambda_{N+1} D +\text{Sym}(Q)+\delta_0 I_m \prec 0.\label{lambdan2}
\end{align} 
thanks to the countability and monotonicity of the eigenvalues of the parabolic operator. The latter implies also that
\begin{align}
-\lambda_{n} D +\text{Sym}(Q)+\delta_0 I_m \prec 0, \quad \forall n\geq N+1. \label{lambdan22}
\end{align}  
We now obtain the following system corresponding to the finite-dimensional part of the eigenspectrum of the parabolic operator:
\begin{align}
\dot{\bar W}(t)=A \bar W(t)+\tilde B v(t),  \label{w^Nopen}
\end{align}
where we denote $\bar W=\text{col}\left\lbrace \bar w_1,\ldots,\bar w_N \right\rbrace \in \mathbb R^{mN}$,
\begin{align} 
    A:=&\text{blkdiag}\{-\lambda_1 D+Q,\ldots,-\lambda_N D+Q \}, \label{Amatrix}
\end{align}
and  $\tilde B\in \mathbb R^{mN\times N}$ is given by
\begin{align}
\tilde B:= \text{col}\left \lbrace B \Psi_1^\top, \ldots, B \Psi_N^\top \right\rbrace =\left ( I_N \otimes B \right )\Psi_{N \times N}.\label{tildeB2}
\end{align}  
{Next, it is more convenient to apply transformation \begin{align} \label{barx} \bar X=\left (I_m \otimes \Psi_{N\times N} \right )X \end{align} for dynamic law \eqref{dynlaw0}. Then, we obtain dynamics }
\begin{align}
\dot{\bar X}(t)=\bar H \bar X(t) -\left (B \otimes \Psi_{N\times N} \right )v(t) \label{dynlaw01}
\end{align}
with 
\begin{align}
\bar H:=&-\bar H_0+Q\otimes I_N;
 \label{barH}\\
\bar H_0:=\text{blkdiag}&\left\lbrace d_1 \Psi_{N\times N}M \Psi_{N\times N}^{-1}, d_2 \Lambda,\ldots, d_m \Lambda \right\rbrace.\notag
\end{align}

 In order to stabilize system \eqref{sys}, we shall select proportional-type actuations $v_j$  to guarantee stabilization of system consisting of \eqref{w^Nopen} and \eqref{dynlaw01}. 
First, we see that by invoking the Hautus lemma, we obtain the following result in conjunction with  the invertibility of $\Psi_{N\times N}$ from Lemma \ref{Lemma_invertible}:
\begin{customlem}{3} \label{lemmacontrollability} 
Let $\psi_j(\cdot)$ be given by \eqref{psi_j} for all $j=1,\ldots,N$. Then, the pair $(A,\tilde B)$ is stabilizable.
\end{customlem}

The abovementioned result guarantees stabilizability of system \eqref{w^Nopen} but not stabilizability of the composite system consisting of \eqref{dynlaw01} and \eqref{w^Nopen}. For this reason, we need to guarantee that matrix $\bar H$ appearing in dynamic control law \eqref{dynlaw01} satisfies a property of the form \begin{align} \text{Sym}(\bar H)\prec -\delta_0 I_{mN}, \ \text{for some } \mu_0 \in \mathbb N \text{ large enough}, \label{SymH}\end{align}
 where $\mu_0$ appears inside $\mu_j$, see \eqref{muj}. The latter is shown to be possible as a result of Assumption \ref{assumption2}.  Based on the above property, in conjunction with  Lemma \ref{lemmacontrollability}, we may choose feedback control law of the form
\begin{align} \label{controlvj}
v_j(t)=K_j \bar W(t), \quad j=1,\ldots,N,
\end{align}
where $K_{j}\in \mathbb R^{1\times mN}$ are controller gains to be found below.  
Property \eqref{SymH} is a result of Assumption \ref{assumption2}, which is rather restrictive contrary to the case of internal stabilization ($\theta=1$), where no particular conditions on the dynamics were imposed. Stabilization of the symmetric part of $H$ is achieved by choice of parameter $\mu_0$ in \eqref{muj}, as it is shown in the following section.

\subsection{Main Boundary Stabilization Result}

We consider here the boundary stabilization of \eqref{sys} by use of dynamic control law \eqref{dynlaw0} and after choice of feedback laws $v_j$. The stability analysis relies on Lyapunov's direct method. 

In order to highlight explicitly the dependence of the feedback control \eqref{controlvj} on dynamic control state $X$ and solution $z$ to \eqref{sys}, we substitute transformations \eqref{transfbarw} and \eqref{transfKar} and we get
\begin{align} 
v_j(t)= &-K_j \Theta  X(t)+K_j \begin{pmatrix} \int_0^L \varphi_1(x) z(t,x) dx\\\vdots\\\int_0^L \varphi_N(x) z(t,x) dx\end{pmatrix}, \label{v_j}
\end{align}
where 
\begin{align} \label{Theta}
\Theta:=\text{blkdiag}\{I_m \otimes \Psi_1,\ldots, I_m \otimes \Psi_N\}.
\end{align}

We are now in a position to establish our main result as a solution to Problem \ref{Problem}(ii) presented in Section \ref{sec:problem}. 
\begin{customthm}{2}\label{proposition} 
Consider parabolic system \eqref{sys} with boundary control ($\theta=0$), boundary conditions satisfying $\gamma_{12}=1-\gamma_{11}$ and $\gamma_{22}=1-\gamma_{21}$ with $\gamma_{11}, \gamma_{21} \in \{0,1\}$, and initial condition $z^0\in H^2\left (0,L;\mathbb R^m\right)$ satisfying $\gamma_{11} z^0(0)+(1-\gamma_{11})\left (z^0\right )^\prime(0)=\gamma_{21}z^0(L)+(1-\gamma_{21})\left (z^0\right )^\prime(L)=0$. Suppose that both Assumption \ref{assumption} on controllability of $(Q,B)$ and Assumption \ref{assumption2} hold true and calculate some $k_Q,\delta_0>0$ satisfying \eqref{stabilizability2}.  Let $N \in \mathbb N$ be subject to \eqref{lambdan2}.  Consider $\mu_j$  and $\psi_j(\cdot)$ given by \eqref{muj}-\eqref{psi_j} with projections $\psi_{j,n}$ given by \eqref{psijn0} and matrix $\Psi_{N\times N}$ as in \eqref{Psi} with inverse explicitly given in Lemma \ref{Lemma_invertible}.
 Assume that there exists $\mu_0 \in \mathbb N$ large enough such that 
\begin{align}
\text{Sym} \left ( \Psi_{N\times N}M\Psi_{N\times N}^{-1}\right )\succ k_Q I_N. \label{ineqM}
\end{align}
with $M=\text{diag}\{\mu_1,\ldots,\mu_N\}.$
Moreover, let $H$ be given by \eqref{Gamma}, $\Theta$ by \eqref{Theta}, and
define
\begin{align}
K=- \Psi_{N\times N}^{-1}\left (I_N\otimes \left (d_1 k_Q B^\top\right )\right).\label{Knproperty01}
\end{align} 
Then, the boundary actuators $u_j(\cdot), j=1,\ldots,N$ subject to dynamic law
\begin{align} \label{dynlawfinal}
\dot X(t)=&\left (H- B \otimes \left (
K\Theta \right )\right ) X(t)\notag\\&+B \otimes
  \left ( d_1 k_Q\Psi_{N\times N}^{-1} \begin{pmatrix}  B^\top \int_0^L \varphi_1(x) z(t,x) dx\\\vdots\\B^\top \int_0^L \varphi_N(x) z(t,x) dx\end{pmatrix} \right ),
\end{align}
(where $X(t)=\text{col}\left\lbrace u(t),r_1(t),\ldots,r_{m-1}(t)\right\rbrace,$ $u(t)=\text{col}\left\lbrace u_1(t),\ldots,u_N(t)\right\rbrace  $) with initial data $X(0)=0$, exponentially stabilize \eqref{sys} with a decay rate $\delta_0$, meaning that the solutions to the closed-loop system satisfy the following inequality: 
\begin{align}
    \Vert z(t,\cdot)\Vert_{L^2\left( 0,L;\mathbb R^m\right)} \leq \ell e^{-\delta_0 t} \Vert z^0(\cdot) \Vert_{L^2\left( 0,L;\mathbb R^m\right)}, \forall t \geq 0 \label{stability01}
\end{align}
with  $\ell>0$.  

Moreover, inequality  \eqref{ineqM} is always feasible for $\mu_0$ large enough. 
\end{customthm} 
\begin{proof}
See Appendix \ref{app:thm2}.
\end{proof}

The abovementioned result illustrates the existence of a  constructive algorithm to stabilize system from the boundary. It mainly relies on the determination of a stabilizing scalar gain $k_Q>0$ corresponding to matrix $Q$ and satisfying \eqref{stabilizability2}, which is independent of the number of unstable modes $N$. It also relies on the determination of parameter $\mu_0$ subject to \eqref{ineqM}. Finally, the inverse of matrix $\Psi_{N\times N}$, given by \eqref{Psi}, is essential to determine the stabilization law, however, we get its explicit formula by Lemma \ref{Lemma_invertible}. Notice also, that \eqref{Knproperty01} here resembles to \eqref{Knproperty0} of Section \ref{sec:internal}, where all $\bar K_j$ there are substituted here by $-d_1 k_Q B^\top$. This allows, when closing the loop of the ODE system \eqref{w^Nopen}, to obtain a matrix with block diagonal elements only, similarly as in \eqref{z^Nsys} of the previous section, where each of the blocks is stabilized by choice of gain $k_Q$. 



\begin{customrem}{4} \label{rem:Hautus}
Note here that for the case of  identical diffusion coefficients, i.e., $d_1=\ldots=d_m$, we might solve the problem by static feedback instead of the dynamic law of Theorem \ref{proposition}. In that case, system is stabilizable in accordance with boundary controllability studies (see \cite{Fernandez-Cara et al (2010)}) as a consequence of the identical diffusion coefficients even if we omit Assumption \ref{assumption2}. Indeed, by performing modal decomposition $z(t,x)=\sum_{n=1}^{+\infty}z_n(t)\varphi_n(x)$, we would obtain the following ODE system for the first $N$ modes: $ \dot Z(t)=A Z(t)+\tilde B u(t)$ with $Z:=\text{col} \left\lbrace z_1,\ldots, z_N\right\rbrace$, $A$ given by \eqref{Amatrix} and $\tilde B= \begin{pmatrix} \tilde B_1 \\ \vdots\\ \tilde B_N  \end{pmatrix} ;\tilde B_n=d_m\left ( (1-\gamma_2)\varphi_n(L)-\gamma_2\varphi_n^\prime(L)\right )B \begin{pmatrix}
1&\cdots&1
\end{pmatrix}.$ Then, since $D=d_mI_m$, the eigenvalues of $A$ are distinct and by the Hautus lemma, system $Z$ is stabilizable, whereas if we had distinct diffusion coefficients, Hautus test would fail. We can, hence, choose proportional controller $u_j(t)=K_j Z(t), j=1\ldots,N$, where $K_j$ are retrieved by a similar inequality as in \eqref{LMIfull0} in Section \ref{sec:internal}. Then, we can follow similar procedure as in the proof of Theorem \ref{theorem1} to show stability of the closed-loop system. However, contrary to the method of Theorem \ref{theorem1} on internal stabilization, we would need here to solve an LMI involving square matrices of dimension $mN$. 
\end{customrem} 

\begin{customrem}{5} 
{
It is worth noting that system \eqref{sys} is a subclass of the general form of controlled systems written abstractly as $\dot y +\mathcal A y =\mathcal Bu$, where $\mathcal A$ has a compact resolvent and a finite number of unstable eigenvectors. Such general classes have been  considered for instance in \cite[Ch. 9]{Munteanu (2019)}. Although constructive methods have been given in these works corresponding to scalar cases, to the best of authors' knowledge, constructive stabilization methods for vector systems have not yet appeared in the literature. The novelty of this work consists in providing completely constructive methods for both internal and boundary stabilization, when the presence of distinct diffusion coefficients complicates the design. This design is based upon modal decomposition combined with Sylvester equations, LMIs and PI controllers.  Recall that our internal stabilization approach provides scalability and relevant independence on the number of unstable modes and it is based on a novel Sylvester-equation approach. For the boundary stabilization approach, we provide a sufficient condition (see Assumption \ref{assumption2}) that leads to the constructive design of a PI controller. For the latter case, there has not appeared a similar approach so far and, to the best of authors' knowledge, only the $2 \times 2$ case has been tackled via backstepping  under stronger sufficient conditions than the ones here \cite{Baccoli et al (2015)}. 
}
\end{customrem} 

\begin{customrem}{6} 
{It would be reasonable to ask why system \eqref{sys} satisfies a cascade form and not a more general form, where $(Q,B)$ would be a controllable pair with no particular structural properties.

(A)  For the internal stabilization case, it turns out that the assumed cascade structure is suitable for the determination of a completely constructive method as in Lemma \ref{lemma} in order to determine transformation \eqref{T_n}. The determination of a similar constructive algorithm for more general pairs $(Q,B)$ would be a very difficult task.  To the best of authors' knowledge, similar transformations as the ones we introduce here have not appeared before. This transformation is subject to easily solvable generalized Sylvester equations that we introduce in this work. Note also that this Sylvester-equation approach is novel in the context of control of PDE theory and it can become a powerful tool not only for this theory but also for (finite-dimensional and large-scale) networked control systems, where simultaneous and scalable stabilization of diagonal systems as in \cite{Dileep et al (2022)} is crucial}. {

(B) For the boundary control case, the chosen cascade form leads to the design of a PI controller after applying transformation \eqref{transfbarw}. For more general cases of the pair $(Q,B)$, finding such a transformation is a difficult task and is left for future research.

Note, however, that if $B$ and $Q$ considered here were both multiplied by a permutation matrix, we would be in a position to follow the same methodology trivially. Such more general pairs would describe alternative systems where the $\bar m$th equation is controlled (with $1\leq \bar m\leq m$), instead of the first equation as in our case. 
}
\end{customrem} 
\section{Simulation}\label{sec:simul}
In the following, we present simulations for both internal and boundary stabilization. These illustrate the results of theorems \ref{theorem1} and \ref{proposition}, respectively. 
\subsection{Internal stabilization example ($\theta=1$)}\label{subsec:exampleinternal}
\begin{figure}
\includegraphics[scale=0.44]{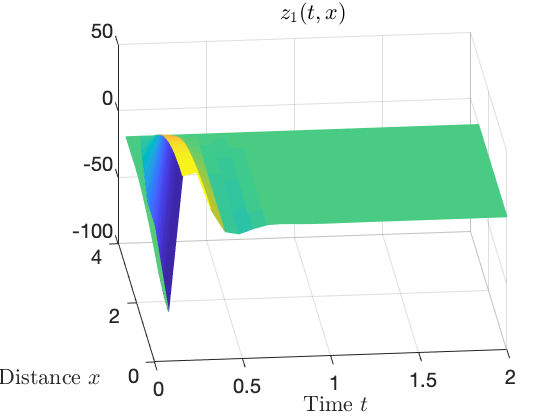}
\caption{Time and space evolution of first state ($\theta=1$)}\label{fig1}
\end{figure}

\begin{figure}
\includegraphics[scale=0.44]{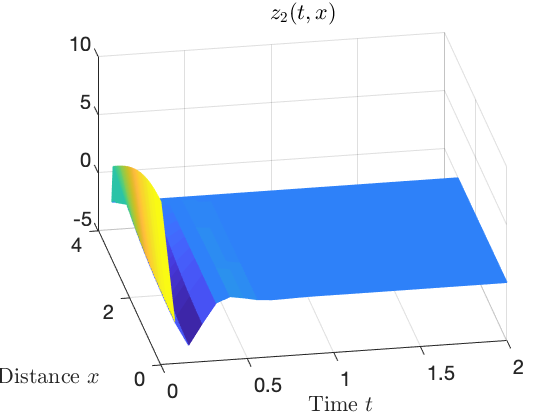}
\caption{Time and space evolution of second state ($\theta=1$) }\label{fig2}
\end{figure}

\begin{figure}
\includegraphics[scale=0.44]{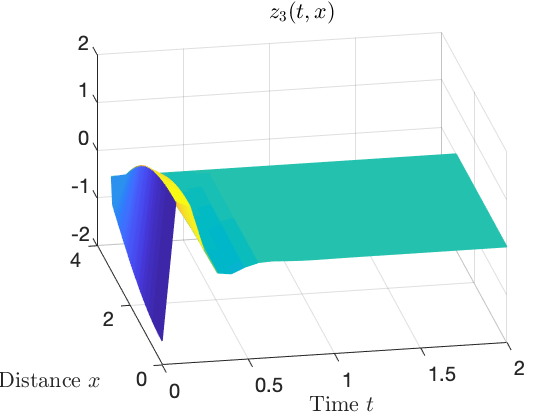}
\caption{Time and space evolution of third state ($\theta=1$)}\label{fig3}
\end{figure}

Let us illustrate the result of Theorem \ref{theorem1} on internal stabilization of \eqref{sys} via an example of $m=3$ PDEs.
Consider $L=\pi$, $\gamma_{11}=0, \gamma_{12}=1, \gamma_{21}=1, \gamma_{22}=0$ meaning that we have Neumann boundary conditions on the left and Dirichlet ones on the right boundary. We choose diffusion matrix and an unstable reaction term given by 
\begin{align}\label{exampledynamics}
D=\text{diag}\{4,3,6 \},\quad   Q=\begin{pmatrix}10&4&8\\ 1&10&2\\0&1&20\end{pmatrix}.
\end{align}
Control is placed internally, hence $\theta=1$.
We retrieve from the Stürm-Liouville problem \eqref{SLprob} the following eigenvalues and eigenfunctions:
\begin{align} \label{lambdaexample}
\lambda_n=\left (n-\frac{1}{2}\right )^2 \pi^2/L^2, \quad \varphi_n(x)=\sqrt 2 \cos(\sqrt{\lambda_n}x). \end{align}
 Let us choose decay rate $\delta=9$.  We select $N=3$ satisfying inequality \eqref{lambdan11}.  Since $m=3$, transformation \eqref{T_n} is of the form \eqref{T_n3} for all $n\leq N$, namely,
$T_n=\begin{pmatrix}
1&3\lambda_n &0\\0&1&0\\0&0&1
\end{pmatrix}.
$
 Shape functions are selected as $b_j(x)=\mathds 1_{[0.1 j, 0.1 j+0.1]}, j=1,\ldots, N,$ in such a way that matrix $\mathcal B_{N\times N}=\begin{pmatrix} 0.1128  &  0.1099 &   0.1042\\
    0.1123  &  0.1059  &  0.0934\\
    0.2226   & 0.1937 &   0.1400\end{pmatrix}$ satisfies Assumption \ref{assumption3} for $N=3$. {We then calculate $K_Q$ by solving LMI \eqref{omega,psi}, which is given by $K_Q=\begin{pmatrix}-67.5&	-3059&-5823\end{pmatrix}$. By invoking \eqref{gains0}, we obtain  $\bar K_1=\begin{pmatrix}-67.3	& -3008 &-5822 \end{pmatrix}, \bar K_2=\begin{pmatrix}-65.3	&-2558	-5809 \end{pmatrix}, $ $\bar K_3=\begin{pmatrix}-61.26&	-1442&	&-5786 \end{pmatrix}.$
We finally calculate the controller gains $K_j$ by using \eqref{Knproperty0}, where $\mathcal B_{N\times N}^{-1}$ $=10^3 \begin{pmatrix} 
4.3906 &  -6.4537    &1.0367\\
   -6.7966 &   9.9474  & -1.5763\\
    2.4255   &-3.5060&    0.5404\end{pmatrix}.$ } Simulations of all three PDE states of the closed-loop system {{} with decay rate $\delta=9$} are shown in figures \ref{fig1}-\ref{fig3} for choice of initial condition $z^0(x)=\left(\begin{pmatrix}\cos x+1\\ 6\cos\frac{x}{2}+3\\-\cos \frac{x}{2}-0.5\end{pmatrix}\right)$.
Furthermore, by using standard LMI solvers, our method illustrated in Theorem \ref{theorem1} via transformation \eqref{T_n} to calculate stabilization gains $K_j$ in \eqref{control0}  is compared with LMI solving resulting from the direct approach \eqref{LMIfull0}. For $N=3$, our indirect approach is approximately $2$ times faster with respect to elapsed time, while for $N=10$ (corresponding to larger $\delta$), it was $200$ times faster than standard LMI.  Note also that for large values of $N$, the non-scalable LMI \eqref{LMIfull0} without our transformation turns out to be computationally hard, while our algorithm to calculate controller gains $K_j$ does not suffer from such limitations. {Note also that our  proposed control method only relies on the inversion of $\mathcal B_{N\times N}$ in \eqref{mathcalB} in order to calculate $K_n$ in \eqref{Knproperty0} and solution of LMI \eqref{omega,psi} of dimension 3. This would not be computationally hard even for large values of $N$, for instance $N=100$, while it would be extremely computationally hard to solve an LMI  \eqref{LMIfull0} of dimension $3N$ (without our transformation).}

\subsection{Boundary stabilization example ($\theta=0$)}

\begin{figure}
\includegraphics[scale=0.44]{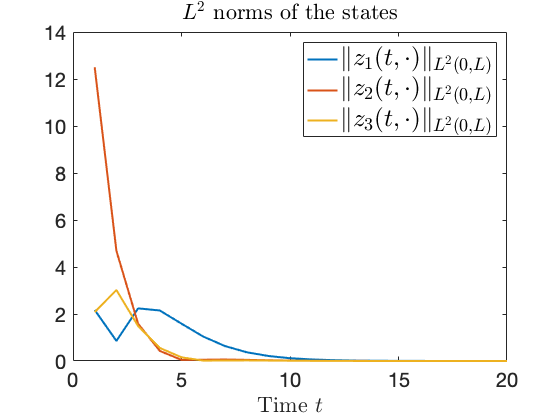}
\caption{Time evolution of the $L^2$ spatial norms of the states ($\theta=0$)}\label{fig4}
\end{figure}

We illustrate here the result of Theorem \ref{proposition} on boundary stabilization ($\theta=0$) of \eqref{sys} via an example of $m=3$ equations. Consider $L=\pi$ and $\gamma_{11}=0, \gamma_{21}=1$ (same boundary conditions as in the previous example of Subsection \ref{subsec:exampleinternal}). Note here that if we select $D$ and $Q$ as in \eqref{exampledynamics}, Assumption \ref{assumption2} would not be satisfied, therefore, we select 
\begin{align}
D=\text{diag}\{4,5,6\}\quad   Q=\begin{pmatrix}10&1&8\\ 1&-10&2\\0&-10&-20\end{pmatrix}. \end{align} 
Eigenvalues and eigenfunctions of the Stürm-Liouville problem \eqref{SLprob} are again as in \eqref{lambdaexample}. We select $k_Q=10$ and $\delta_0=9$ satisfying Assumption \ref{assumption2}. We select $N=3$ for which inequality \eqref{lambdan2} is satisfied. Also, select $\mu_0=5$ in such a way that \eqref{ineqM} is satisfied. Then, $\mu_j$ satisfy $\sqrt{\mu_j}=j\frac{\pi}{L}+2\mu_0\frac{\pi}{L}$ (see \eqref{muj}). Functions $\psi_j(\cdot)$ (see \eqref{psi_j}) are given by $\psi_j(x)= (-1)^j \cos\left ( \sqrt{\mu_j} x\right )$. Matrices $\Psi_{N\times N}$ and its inverse are given by (see \eqref{Psi} and \eqref{invPsi})
$\Psi_{N\times N}= \begin{pmatrix} -0.0064 &  -0.0054 &  -0.0047\\
    0.0194  &  0.0165  &  0.0143\\
   -0.0334  & -0.0283  & -0.0243\end{pmatrix},$ $\Psi_{N\times N}^{-1}=10^7\begin{pmatrix}
   -0.3987  &-0.1945 &  -0.0370\\
    1.0829    &0.5268  &  0.0996\\
   -0.7110   &-0.3452   &-0.0650\end{pmatrix} .$ We then apply Theorem \ref{proposition} by considering dynamic law as in \eqref{dynlawfinal} with $K=- 40\Psi_{N\times N}^{-1}\left (I_N\otimes  B^\top\right)$. In Figure \ref{fig4}, we see the evolution of the $L^2$ spatial norms of all three system states {with decay rate $\delta_0=9$ } for choice of initial condition $z^0(x)=\left(\begin{pmatrix}\cos x+1\\ 6\cos\frac{x}{2}+3\\-\cos \frac{x}{2}-0.5\end{pmatrix}\right)$. 

\section{Conclusion} \label{conclusion}
The problem of both internal and boundary stabilization of an underactuated parabolic system in a cascade form and in the presence of distinct diffusion coefficients was considered. For the internal stabilization problem, after performing modal decomposition, the problem was reduced to just the stabilization of the reaction term avoiding in that way a direct stabilization of the whole system of ODEs corresponding to the comparatively unstable modes, which might have arbitrarily large dimension. An easily calculable state transformation of dimension equal to the number of coupled PDEs as a solution to a generalized Sylvester equation was introduced in order to solve this stabilization reduction problem. For the boundary stabilization problem, we used a dynamic extension technique in such a way that the control components are placed internally in the PDEs. Again, the stabilization problem was reduced to just the stabilization of the reaction term.

In our future works, the present approach will be applied to observer-based control and extended to nonlinear systems.

\appendices
\section{Proof of Theorem \ref{theorem1}} \label{appendix1}

Below, we prove Theorem \ref{theorem1} on internal stabilization of Section \ref{sec:internal}. Note first that transformation $T_n$ appearing in stabilization law is calculated via the constructive Algorithm \ref{algorithm} coming from Lemma \ref{lemma}. To see how $T_n$ maps \eqref{zn_dot0} to target system \eqref{target_sys}  via control \eqref{control0}, let us apply it to \eqref{zn_dot0}. Therefore, we obtain
\begin{align}
    \dot y_n(t)=&\left (-\lambda_n T_n D T_n^{-1 } +T_n Q T_n^{-1}+B K_n T_n^{-1}\right )y_n(t),
\end{align}
for all $n=1,\ldots,N.$
Comparing the above system with target system \eqref{target_sys}, the following equations must be satisfied for all $n=1,\ldots,N$:
\begin{align}
    \left ( Q-\lambda_n d_m I_m\right )T_n +T_n\left ( \lambda_n D-Q\right )+B G_n T_n=0. \notag
\end{align}
Substituting  \eqref{T_n} in the previous equation, this is written as
\begin{align}
   &\left ( Q-\lambda_n d_m I_m\right )\left( \sum_{i=1}^{\bar \sigma} \bar T_i \lambda_n^i\right )+\left( \sum_{i=1}^{\bar \sigma} \bar T_i \lambda_n^i\right )\left ( \lambda_n D-Q\right )\notag\\&+\left ( D-d_mI_m\right )\lambda_n+B G_n T_n =0, \quad n=1,\ldots,N. \label{sylvesterwithGn}
\end{align}
After injecting expression for $G_n$, \eqref{sylvesterwithGn} is written as 
\begin{align}
    &\left( I_m-BB^\top\right ) \left ( \left ( Q-\lambda_n d_m I_m\right )\left( \sum_{i=1}^{\bar \sigma} \bar T_i \lambda_n^i\right )\right.\notag\\&\left. +\left( \sum_{i=1}^{\bar \sigma} \bar T_i \lambda_n^i\right )\left ( \lambda_n D-Q\right )+\left ( D-d_mI_m\right )\lambda_n\right )=0. \label{firstSYlvester}
\end{align}
Then, eliminating all the coefficients of $\lambda_n^i$ in \eqref{firstSYlvester} for all $i$ in $\{1,\ldots,\bar \sigma\}$, we obtain \eqref{Sylvester}, which is assumed to hold true for all $\bar T_i$, $i\in \{1,\ldots,\bar \sigma\}$. Therefore, \eqref{Sylvester} guarantees that, via transformation $T_n$, we obtain target system \eqref{target_sys}.

At this point, let us remark that for given initial condition $z^0$ in $H^1\left (0,L;\mathbb R^m\right)$ satisfying compatibility conditions for \eqref{sys}, unique existence of classical solutions to system \eqref{sys} with nonlocal feedback control \eqref{control0} ($\theta=1$), namely $z\in C^1\left ([0,+\infty);L^2\left (0,L;\mathbb R^m\right) \right )$, follows from simple argument such as the Lumer-Philipps theorem, see for example (\cite{Pazy (1983)}, Corollary 4.4, Chapter 1.  

Let us now prove $L^2$ stability of the closed-loop system applying direct Lyapunov method (see for example \cite{Katz and Fridman (2020)}). First, observe that by injecting gains \eqref{gains0}, $Y:=\text{col}\{ y_1,\ldots,y_N\} \in \mathbb R^{mN}$ satisfies dynamic law \begin{align} 
    \dot Y(t)=R Y(t), \label{y^Nclosed}
\end{align}
where $R:=\text{blkdiag}\{R_1,\ldots,R_N\}$ with $R_n:=-\lambda_n d_mI_m+ Q+B K_Q, n=1,\ldots,N.$
Now, by the fact that $(Q,B)$ is controllable, we can stabilize matrix $Q$, in such a way that we can find $0 \prec P \in \mathbb R^{m \times m}$, $K_Q\in \mathbb R^{1\times m}$, and $\bar \rho>0$ such that
\begin{align}
\scalemath{0.95}{\begin{pmatrix} \text{Sym} \left ( P(Q+B K_Q) \right )+\left(\delta-\lambda_1 d_m \right ) P &I_m\\I_m&-\bar \rho I_m \end{pmatrix}\prec 0, }\label{LMIP}
\end{align}
which is written in the design LMI form
\begin{align} 
\scalemath{1}{\begin{pmatrix} \text{Sym}\left (QP_{-1}+BJ\right )+\left(\delta-\lambda_1 d_m \right ) P_{-1} &I_m\\I_m&-\bar \rho I_m\end{pmatrix} \prec 0,} \notag
\end{align} 
where $P_{-1}=P^{-1}$ and $J=K_Q P_{-1}$. The latter implies also feasibility of LMI \eqref{omega,psi}.
Next, by virtue of \eqref{lambdan11}, we can always find $\rho>0$ such that the following LMI is satisfied:
\begin{align}
\begin{pmatrix}
-\lambda_{N+1} D +\text{Sym}(Q)+\delta I_m & \frac{1}{\sqrt 2}I_m\\\frac{1}{\sqrt 2}I_m&-\rho I_m
\end{pmatrix}\prec 0.  \label{Omega}
\end{align}
To prove stability, defining first $y(t,x)=\sum_{n=1}^{+\infty}\varphi_n(x)y_n(t)$, consider Lyapunov functional $\mathcal V:L^2\left (0, L;\mathbb R^m\right)\to \mathbb R$
\begin{align}
\mathcal V[y]=\frac{1}{2}  Y^\top \bar P Y+\frac{\eta}{2} \sum_{n=N+1}^{+\infty} \vert y_n\vert^2,
\end{align} \begin{align}\label{rho0}\text{where }\eta:=\frac{2}{\rho \bar \rho\beta \vert K \vert^2};  \beta:=\max_{n=1,\ldots,N} \vert T_n^{-1}\vert^2 \sum_{j=1}^{N}\Vert b_j(\cdot) \Vert_{L^2 \left (0,L\right )}^2\end{align} with $\rho>0$ satisfying \eqref{Omega}, $\bar \rho>0$ satisfying \eqref{LMIP} and $K$ given by \eqref{Knproperty0} (where $\bar K_n$ are given by \eqref{gains0}). Also, $\bar P:=I_N \otimes P$.  By invoking boundedness of $\mathcal T, \mathcal T^{-1}$ in $\ell^2\left (\mathbb N;\mathbb R^{m}\right)$, the fact that $y_n=T_n z_n,$ the Cauchy-Schwarz inequality, and Parseval's identity, we get $\underline c \Vert z(t,\cdot)\Vert^2_{L^2 \left (0,L;\mathbb R^m\right )}=\underline c \sum_{n=1}^{+\infty}\vert z_n(t)\vert^2\leq \sum_{n=1}^{+\infty} \vert y_n(t)\vert^2 \leq  \bar c \Vert z(t,\cdot)\Vert^2_{L^2 \left (0,L;\mathbb R^m\right )}, $
where $\underline c:=\frac{1}{\max_{n\in \mathbb N}\vert T_n^{-1}\vert^2}$ and $\bar c:=\max_{n \in \mathbb N}\vert T_n \vert^2.$
By continuous differentiability of solutions with respect to $t$ for all $t\geq 0$, we are in a position to define $V(t):=\mathcal V[y](t)$ for all $t\geq 0$ and we may take its time-derivative $\dot V(t)$ along the solutions of target system \eqref{target_sys}. By use of the previous inequality, we obtain for $V(t)$
\begin{align}
\frac{\underline c}{2}\min \left (\lambda_{\min}(P),\eta\right )\Vert z(t,\cdot) \Vert_{L^2 \left (0,L;\mathbb R^m\right )}^2\leq V(t) \notag\\\leq \frac{\bar c}{2}\max \left (\lambda_{\max}(P),\eta\right )\Vert z(t,\cdot) \Vert_{L^2 \left (0,L;\mathbb R^m\right )}^2. \label{Vineq}
\end{align}
Its derivative satisfies
\begin{align}
\dot V(t) =& Y^\top(t) \text{Sym}(\bar P R)Y(t)\notag\\&+\eta\sum_{n=N+1}^{+\infty}y_n^\top(t) \left( -\lambda_n D + \text{Sym}\left ( Q\right )\right) y_n(t)\notag\\&+\eta \sum_{n=N+1}^{+\infty}y_n^\top(t) B \sum_{j=1}^N b_{j,n} K_j Z(t). \label{Lyapder}
\end{align}
By the Cauchy-Schwarz inequality and Parseval's identity, last term of \eqref{Lyapder} is bounded as follows:
\begin{align*}
&\eta \sum_{n=N+1}^{+\infty}y_n^\top(t) B \sum_{j=1}^N b_{j,n} K_j Z(t) \leq \eta \frac{1}{2\rho}  \sum_{n=N+1}^{+\infty} \vert y_n(t)\vert^2 \notag\\&+\eta\frac{\rho}{2}  \vert Z(t)\vert^2 \vert K\vert^2 \sum_{n=N+1}^{+\infty}\vert \mathcal B_n^\top \vert^2\notag\\&\leq \eta\frac{1}{2\rho} \sum_{n=N+1}^{+\infty} \vert y_n(t)\vert^2 +\eta \frac{\rho}{2}\beta \vert K\vert^2 \vert Y(t)\vert^2
\end{align*}
where $\rho>0$ satisfies \eqref{Omega}  and $\beta$ is given by \eqref{rho0}. After substituting expression \eqref{rho0} for $\eta$, \eqref{Lyapder} is bounded as
{\small{\begin{align}
\dot V(t)\leq -2 \delta V(t) +Y^\top(t) \Gamma Y(t)+\eta\sum_{n=N+1}^{+\infty} y_n^\top(t) \Omega_n y_n(t) \label{Lyapder20}
\end{align} } }
 \begin{align*}\text{where }\Gamma:=\text{blkdiag}\{\text{Sym}(P R_1)+\left ( \frac{1}{\bar \rho}+\delta P \right )I_{m},\ldots,\\\text{Sym}(P R_N)+\left ( \frac{1}{\bar \rho}+\delta P \right )I_{m}\},\end{align*} $$\Omega_n:=-\lambda_n D+\text{Sym}(Q)+\left (\frac{1}{2\rho}+\delta \right )I_m.$$
Monotonicity of the eigenvalues, in conjunction with \eqref{LMIP} and \eqref{Omega}, implies $\Gamma<0$ and $ \Omega_n<0, \quad \forall n \geq N+1,$ respectively. Thus, \eqref{Lyapder20} in conjunction with \eqref{Vineq} readily yields to a stability inequality of the form \eqref{stability0}. 

The proof of Theorem \ref{theorem1} is complete.  $\blacksquare$

\section{Proof of Theorem \ref{proposition}}\label{app:thm2}
We prove here Theorem \ref{proposition} in Section \ref{sec:boundary} on boundary stabilization.

We invoke first existence-uniqueness of solutions to the closed loop system \eqref{sys}, \eqref{dynlaw0} with $\theta=0$ by easily adapting a result given in \cite{Karafyllis (2021)} for the scalar case to our vector case (proof of Theorem 2.2 therein). More precisely, for any given initial condition $z^0\in H^2\left (0,L;\mathbb R^m\right)$ satisfying $\gamma_{11} z^0(0)+(1-\gamma_{11})\left (z^0\right )^\prime(0)=\gamma_{21}z^0(L)+(1-\gamma_{21})\left (z^0\right )^\prime(L)=0$ (implying by \eqref{transfKar} that $w^0(\cdot):=w(0,\cdot) \in H^2\left (0,L;\mathbb R^m\right)$ satisfying $\gamma_{11} w^0(0)+(1-\gamma_{11})\left (w^0\right )^\prime(0)=\gamma_{21}w^0(L)+(1-\gamma_{21})\left (w^0\right )^\prime(L)=0$) and input initial conditions $u(0)=0$, there exists a unique solution $w\in C^0\left ( [0,+\infty)\times[0,L];\mathbb R^m\right )\cap C^1\left ( (0,+\infty)\times[0,L];\mathbb R^m\right )$ with $w(t,\cdot)\in C^2\left ([0,L];\mathbb R^m \right )$ of the closed loop system \eqref{wsystem}, \eqref{dynlawfinal} implying also unique existence of $z$ in the same function spaces due to \eqref{transfKar}. Simultaneously, we get $u\in C^1\left ( [0,+\infty);\mathbb R^N\right )$.

Next, notice that dynamic law \eqref{dynlawfinal} is directly deduced by \eqref{dynlaw0} after substituting expression for $v_j$ in \eqref{v_j} and also gains \eqref{Knproperty01}. We also see that inequality \eqref{ineqM} is feasible for choice of $\mu_0$ large enough. Indeed, $\mu_j$ given by \eqref{muj} are written in the form $\mu_j=\mu_0^2\left ( 4 \frac{\pi^2}{L^2}+4\frac{\pi}{L}\sqrt{\bar \mu_j}\frac{1}{\mu_0}+\frac{\bar \mu_j}{\mu_0^2}\right ), j=1,\ldots,N$. Then, recalling that $M:=\text{diag}\{\mu_1,\ldots,\mu_N\},$ it is easy to see that whenever $\mu_0 \to +\infty$, we obtain $\text{Sym} \left ( \Psi_{N\times N}M\Psi_{N\times N}^{-1}\right )=O(\mu_0^2)I_N,$
which yields feasibility of \eqref{ineqM}. 

In the next step, let us prove $L^2$ stability of the closed-loop system \eqref{barwn1}, \eqref{dynlaw01} by applying direct Lyapunov method. First, observe that by injecting control law \eqref{controlvj} and gains \eqref{Knproperty01} in \eqref{w^Nopen} and \eqref{dynlaw01}, $\bar W$ and $\bar X$ satisfy the following dynamics:
\begin{subequations}\label{finalsystem}
 \begin{align}
\begin{aligned} 
    \dot{\bar W}(t)=&  \Pi_{1} \bar W(t),\\
\dot{\bar X}(t)=&   \bar H\bar X(t)+ \Pi_{2} \bar W(t)
    \end{aligned}
\end{align}
where 
\begin{align}
 \Pi_{1}:=&\text{blkdiag}\{Q-D\text{diag}\{k_Q, \lambda_n,\ldots,\lambda_n\}\}_{n=1}^N,\notag\\
\Pi_{2}:=&d_1k_Q B \otimes \left ( I_N \otimes B^\top\right ),
\end{align}
and $\bar H$ is given by \eqref{barH}.
In addition, by \eqref{barwn1}, we get the following dynamics for all $n\geq N+1:$
\begin{align}
\dot{\bar w}_n(t)=\left (-\lambda_n D +Q \right )\bar w_n(t)+B\sum_{j=1}^N \psi_{j,n} K_j \bar W(t),
\end{align}
\end{subequations}
Next, see that by by virtue of \eqref{lambdan2}, we can always find a $\rho>0$ such that the following LMI is satisfied:
\begin{align}
\begin{pmatrix}
-\lambda_{N+1} D +\text{Sym}(Q)+\delta_0 I_m & \frac{1}{\sqrt 2}I_m\\\frac{1}{\sqrt 2}I_m&-\rho I_m
\end{pmatrix}\prec 0.  \label{Omega2}
\end{align}
Now, notice that by Assymption \ref{assumption2}, \begin{align}&\text{Sym} \{ \Pi_{1}\}=\text{blkdiag}\{\text{Sym}(Q) -D\text{diag}\{k_Q, \lambda_n,\ldots,\lambda_n\}\}_{n=1}^N\notag\\&\preceq I_N \otimes \left (\text{Sym}(Q) -D\text{diag}\{k_Q, \lambda_1,\ldots,\lambda_1\}\right )\preceq -\delta_0 I_{mN},\notag\end{align}
from which we can always find a $\bar \rho>0$ such that 
 \begin{align}&\text{Sym}( \Pi_{1})+\frac{1}{\bar \rho}I_{mN}\prec -\delta_0 I_{mN}.\label{Pi11ineq2}\end{align}
Also, by invoking \eqref{ineqM} and by virtue of Assumption \ref{assumption2}, we get
\begin{align}
&\text{Sym}(\bar H) \notag\\&=-\text{blkdiag}\left\lbrace d_1  \text{Sym}\left (\Psi_{N\times N}M  \Psi_{N\times N}^{-1}\right ), d_2 \Lambda ,\ldots,d_m\Lambda \right\rbrace\notag\\&+\text{Sym}(Q) \otimes I_N \notag\\&\preceq -\text{blkdiag}\left\lbrace d_1  k_Q I_N, d_2 \lambda_1 I_N,\ldots,d_m\lambda_1 I_N \right\rbrace\notag\\&+\text{Sym}(Q) \otimes I_N \preceq -\delta_0 I_{mN}, \notag
\end{align}
from which we can always find a $\bar{\bar \rho}>0$ such that 
\begin{align}
\text{Sym}(\bar H) +\frac{1}{\bar{\bar \rho}}I_{mN}\prec -\delta_0 I_{mN}. \label{Pi22ineq2}
\end{align}
The above is a desired property as it was already mentioned in \eqref{SymH}.
To prove stability, defining $\bar w(t,x)=\sum_{n=1}^{+\infty} \phi_n(x)\bar w_n(t)$, consider Lyapunov functional $\mathcal V:L^2\left (0, L;\mathbb R^m\right)\times \mathbb R^{mN} \to \mathbb R$
{\small{\begin{align}
\mathcal V[\bar w, \bar X]=&\frac{1}{2} \vert \bar W\vert^2+\frac{\eta_1}{2} \sum_{n=N+1}^{+\infty} \vert \bar w_n\vert^2+\frac{\eta_2}{2} \vert \bar X\vert^2, \text{ where }\\
\eta_1:=&\frac{1}{\bar \rho \rho \sum_{j=1}^N \Vert \psi_j(\cdot)\Vert_{L^2(0,L)}^2 \vert K\vert^2}, 
\eta_2:=\frac{1}{\bar{\bar{\rho}} \rho d_1^2 k_Q^2}. \label{etas}
\end{align}}}
By  use of transformations \eqref{transfbarw}, \eqref{barx} and Parseval's identity, we obtain
$\Vert w(t,\cdot)\Vert^2_{L^2 \left (0,L;\mathbb R^m\right )}+\underline c \vert X(t)\vert^2$ $= \sum_{n=1}^{+\infty}\vert w_n(t)\vert^2+\underline c \vert X(t)\vert^2\leq \sum_{n=1}^{+\infty} \vert \bar w_n(t)\vert^2 +\vert \bar X(t)\vert^2$ $\leq  \Vert w(t,\cdot)\Vert^2_{L^2 \left (0,L;\mathbb R^m\right )}+\bar c\vert X(t)\vert^2,
$
where $\underline c:=\frac{1}{\vert \Psi_{N\times N}^{-1}\vert^2 }+ \sum_{n=1}^{+\infty}\vert \Psi_n\vert^2=\frac{1}{\vert \Psi_{N\times N}^{-1}\vert^2 }+ \sum_{j=1}^{N}\Vert \psi_j(\cdot)\Vert_{L^2(0,L)}^2$ and $\bar c:=\vert \Psi_{N\times N}\vert^2+\sum_{j=1}^{N}\Vert \psi_j(\cdot)\Vert_{L^2(0,L)}^2$. 
 By continuous differentiability of solutions with respect to $t$ for all $t\geq 0$, we are in a position to define $V(t):=\mathcal V[\bar w,\bar X](t)$ for all $t\geq 0$ and we may take its time-derivative $\dot V(t)$ along the solutions of system \eqref{finalsystem}. By use of the previous inequality, we obtain for $V(t)$
\begin{align}
&\underline C_1 \Vert w(t,\cdot) \Vert_{L^2 \left (0,L;\mathbb R^m\right )}^2+\underline C_2 \vert X(t)\vert^2 \leq V(t) \notag\\&\leq \bar C_1 \Vert w(t,\cdot) \Vert_{L^2 \left (0,L;\mathbb R^m\right )}^2+\bar C_2 \vert X(t)\vert^2, \label{Vineq2}
\end{align}
where $\underline C_1:=\frac{1}{2}\min\{1,\eta_1\}, \underline C_2:=\eta_2\frac{ \underline c }{2} $ and $\bar C_1:=\frac{1}{2}\max\{1,\eta_1\}, \bar C_2:=\eta_2\frac{ \bar c }{2}.$ 
Differentiating $V$ along the solutions of \eqref{finalsystem}, we obtain
{\small{\begin{align}
&\dot V(t) =  \bar W^\top(t) \text{Sym}(\Pi_{1} )\bar W(t)\notag\\&+\eta_1 \sum_{n=N+1}^{+\infty}\bar w_n^\top(t) \left( -\lambda_n D + \text{Sym}\left ( Q\right )\right) \bar w_n(t)\notag\\&+\frac{\eta_1}{2} \sum_{n=N+1}^{+\infty}\bar w_n^\top(t) B \psi_{j,n} K_j \bar W(t)\notag\\&+\frac{\eta_1}{2}\bar W^\top(t) \sum_{n=N+1}^{+\infty} K_j^\top \psi_{j,n} B^\top \bar w_n(t)\notag\\&+\eta_2 \bar X^\top(t) \text{Sym}\left (\bar H\right ) \bar X(t)+\frac{\eta_2}{2} \bar X^\top(t)  \Pi_{2}\bar W(t)\notag\\&+\frac{\eta_2}{2}\bar W^\top (t){\Pi}_{2}^\top\bar X(t). \label{Lyapder2}
\end{align}}}
We apply next Young's inequality for the cross terms in the Lyapunov derivative as follows:
{\small{\begin{align}
&\dot V(t) \leq  \bar W^\top(t) \text{Sym}(\Pi_{1} )\bar W(t)\notag\\& +\eta_1 \sum_{n=N+1}^{+\infty}\bar w_n^\top(t) \left( -\lambda_n D + \text{Sym}\left ( Q\right )\right) \bar w_n(t)\notag\\&+\eta_1 \frac{1}{2\rho} \sum_{n=1}^{+\infty} \vert \bar w_n\vert^2+\eta_1 \frac{\rho}{2}\sum_{j=1}^N \Vert \psi_j(\cdot)\Vert_{L^2(0,L)}^2 \vert K\vert^2 \vert \bar W(t)\vert^2\notag\\&+\eta_2 \bar X^\top(t) \text{Sym}\left (\bar H\right ) \bar X(t)+{\eta_2}\frac{1}{2\bar{\bar \rho}}\vert \bar X(t)\vert^2+\eta_2\frac{\bar{\bar \rho}}{2} d_1^2 k_Q^2  \vert \bar W(t)\vert^2.\notag
\end{align} }}
By substituting \eqref{etas} and by use of \eqref{Omega2}, \eqref{Pi11ineq2}, and \eqref{Pi22ineq2}, we obtain $\dot V(t) \leq  -\delta_0\vert \bar W(t) \vert^2 $ $-\eta_1 \delta_0\sum_{n=N+1}^{+\infty} \vert \bar w_n(t)\vert^2-\eta_2 \delta_0 \vert \bar X(t)\vert^2=-2\delta_0 V(t), \quad \forall t \geq 0. $
Finally, by combining the previous inequality with \eqref{Vineq2}, the fact that $X(0)=0$,  and the fact that from transformation \eqref{transfKar}, we have $\Vert z(t,\cdot)\Vert^2_{L^2\left ( 0,L;\mathbb R^m\right )} \leq \Vert w(t,\cdot)\Vert^2_{L^2\left ( 0,L;\mathbb R^m\right )} $ $+ \sum_{j=1}^N \Vert \psi_j(\cdot)\Vert^2_{L^2\left ( 0,L\right )} \vert X(t)\vert^2,$
we readily obtain \eqref{stability01}. $\blacksquare$

\begin{IEEEbiography}[{\includegraphics[width=1.10in,height=1.10in,clip]{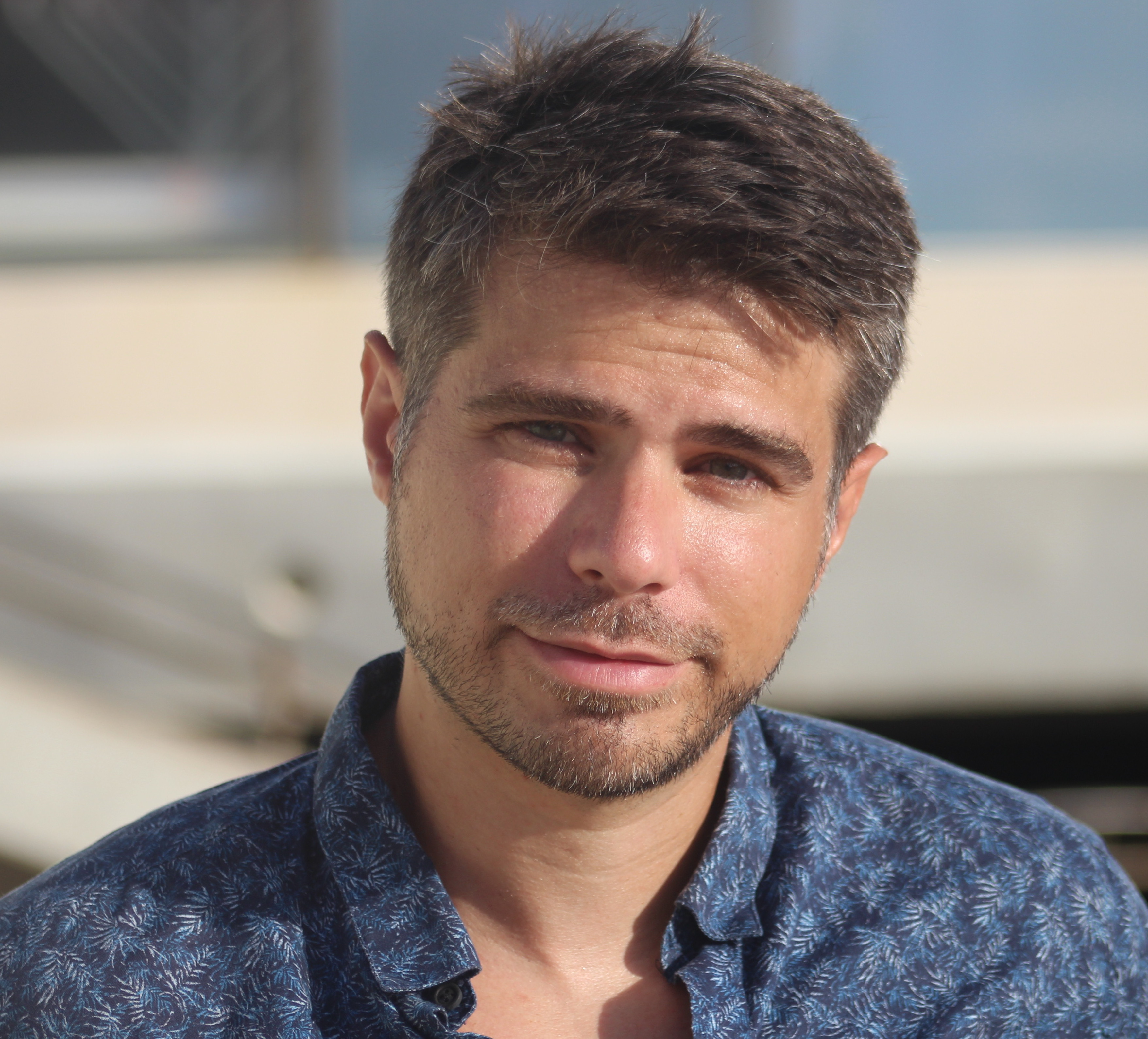}}]
{Constantinos Kitsos} holds a Diploma in Electrical and Computer Engineering and in 2016 he received an M.Sc. in Applied Mathematics, both from the National Technical University of Athens, Greece. In 2020, he received a Ph.D. degree in Automatic Control from Université Grenoble Alpes (GIPSA-lab), France. Since then, he has been affiliated as a Postdoctoral Researcher with the Department of Electrical Engineering of Tel-Aviv University, Israel and with with the Laboratory for Analysis and Architecture of Systems of the French National Center for Scientific Research (LAAS-CNRS), Toulouse, France. His research interests include nonlinear observers and  control of PDEs.
\end{IEEEbiography}
\vskip -2\baselineskip plus -1fil
\begin{IEEEbiography}[{\includegraphics[width=1.13in,height=1.10in,clip]{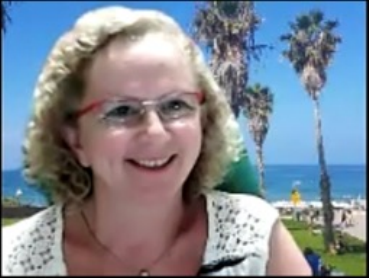}}]
{Emilia Fridman} received the M.Sc. and the Ph.D. degrees in USSR in mathematics. Since 1993 she has been at Tel Aviv University, where she is currently Professor of Electrical Engineering-Systems. She has held visiting positions at the Weierstrass Institute for Applied Analysis and Stochastics in Berlin (Germany), INRIA in Rocquencourt (France), Ecole Centrale de Lille (France), Leicester University (UK), Kent University (UK), CINVESTAV (Mexico), Zhejiang University (China), St. Petersburg IPM (Russia), Melbourne University (Australia), Supelec (France), KTH (Sweden).

Her research interests include time-delay systems, networked control systems, distributed parameter systems, robust control, singular perturbations and nonlinear control. She has published two monographs and more than 200 articles in international scientific journals. She serves/served as Associate Editor in Automatica, SIAM Journal on Control and Optimization and IMA Journal of Mathematical Control and Information. In 2014 she was recognized as a Highly Cited Researcher by Thomson ISI. Since 2018, she has been the incumbent for Chana and Heinrich Manderman Chair on System Control at Tel Aviv University. She is IEEE Fellow since 2019. In 2021 she was recipient of IFAC Delay Systems Life Time Achievement Award and of Kadar Award for outstanding research in Tel Aviv University. She was a member of the IFAC Council. She is currently IEEE CSS Distinguished Lecturer. In 2023 her monograph ``Introduction to Time-Delay Systems: Analysis and Control" (Birkhauser, 2014) was the winner of IFAC Harold Chestnut Control Engineering Textbook Prize.
\end{IEEEbiography}

\end{document}